\documentclass[a4paper]{amsart}

\usepackage[all]{xy}
\usepackage{latexsym}
\usepackage{amssymb}
\usepackage{amsmath}
\usepackage{amsthm}
\usepackage{amscd}
\usepackage{fancyhdr}
\usepackage[dvips]{graphicx}
\usepackage{a4wide}
\usepackage{fancybox}
\usepackage[mathscr]{eucal}
\usepackage{mathtools}
\usepackage{pdfsync}
\usepackage{color}
\usepackage{colordvi}
\definecolor{Chocolat}{rgb}{0.36, 0.2, 0.09}
\definecolor{BleuTresFonce}{rgb}{0.215, 0.215, 0.36}
\usepackage[colorlinks,final,backref=page,hyperindex]{hyperref}
\hypersetup{citecolor=BleuTresFonce, linkcolor=Chocolat}

\theoremstyle{plain}
\newtheorem{pro}{Proposition}[section]
\newtheorem{thm}[pro]{Theorem}

\newtheorem{lem}[pro]{Lemma}

\newtheorem{cor}[pro]{Corollary}

\newtheorem*{thmA}{Theorem A}
\newtheorem*{thmB}{Theorem B}

\theoremstyle{definition}
\newtheorem{defi}[pro]{Definition}

\newtheorem*{oldefi}{Definition}
\theoremstyle{remark}
\newtheorem*{Rq}{\sc Remark}

\newtheorem*{Ex}{\sc Example}

\newcommand{\bialg}{H}
\newcommand{\prdalg}{\mu}
\newcommand{\coprdalg}{\varDelta}
\newcommand{\unalg}{\mathit{u}}
\newcommand{\counalg}{\varepsilon}
\newcommand{\SHMod}{\gsym\text{-}\bialg\textsf{-Mod }}

\newcommand{\dgSHMod}{\textsf{dg} \ \gsym\text{-}\bialg\textsf{-Mod} }
\newcommand{\HMod}{\bialg\textsf{-Mod }}
\newcommand{\dgHMod}{\textsf{dg }\bialg\textsf{-Mod }}

\newcommand{\HOpd}{\mathsf{Op}_{\bialg}}

\newcommand{\HCoopd}{\mathsf{Coop}_{\bialg}}

\newcommand{\Sw}[2]{#1_{(#2)}}
\newcommand{\gsym}{\mathbb{S}}
\newcommand{\SH}{$\gsym$-$\bialg$}
\newcommand{\opd}{\mathcal{P}}
\newcommand{\comp}{\gamma}
\newcommand{\unit}{\eta}
\newcommand{\coopd}{\mathcal{C}}
\newcommand{\decomp}{\Delta}
\newcommand{\coun}{\eta}
\newcommand{\anti}{S}
\newcommand{\alg}{A}
\newcommand{\KK}{\mathbb{K}}
\newcommand{\NN}{\mathbb{N}}

\newcommand{\Ical}{\mathcal{I}}

\newcommand{\freeop}{\mathcal{T}}
\newcommand{\End}{\mathrm{End}}
\newcommand{\Hom}{\mathrm{Hom}}
\newcommand{\Id}{\mathrm{Id}}
\newcommand{\id}{\mathrm{id}}
\newcommand{\infprd}[2]{#1 \circ_{(1)} #2}
\newcommand{\twm}{\alpha}
\newcommand{\Tw}{\mathrm{Tw}}

\newcommand{\Barc}{\mathrm{B}}
\newcommand{\Cobc}{\Omega}
\newcommand{\gensp}{E}
\newcommand{\rel}{R}

\newcommand{\hadprd}{\underset{\textrm{H}}{\otimes}}
\newcommand{\BV}{\mathcal{BV}}

\newcommand{\Gerst}{\mathcal{G}}
\newcommand{\prd}[2]{#1\bullet #2}
\newcommand{\bracket}[2]{\langle #1  ;  #2 \rangle}
\newcommand{\unop} [1]{\Delta #1}
\newcommand{\unopalg}{\delta}
\newcommand{\algnbduaux}{D}
\newcommand{\Com}{\mathcal{C}om}
\newcommand{\Lie}{\mathcal{L}ie}
\newcommand{\suspopd}{\mathcal{S}}
\newcommand{\underlying}[1]{\underline{#1}}

\makeatletter
\let\original@addcontentsline\addcontentsline
\newcommand{\dummy@addcontentsline}[3]{}
\newcommand{\DeactivateToc}{\let\addcontentsline\dummy@addcontentsline}
\newcommand{\ActivateToc}{\let\addcontentsline\original@addcontentsline}
\makeatother

\newenvironment{proo}{\begin{trivlist} \item{\sc {Proof.}}}
  { \begin{flushright}
 $\square$             \end{flushright} \end{trivlist}}

\long\def\symbolfootnote[#1]#2{\begingroup%
\def\thefootnote{\fnsymbol{footnote}}\footnote[#1]{#2}\endgroup}

\begin{document}

\title{Koszul duality theory for operads over Hopf algebras}
\author{Olivia BELLIER}
\address{Laboratoire J.-A. Dieudonn\'e \\ Universit\'e de Nice - Sophia Antipolis \\ Parc Valrose \\ 06108 NICE Cedex 02 \\ FRANCE}
\curraddr{Institut de Math\'ematiques de Toulouse \\ Universit\'e Paul Sabatier \\ 118, route de Narbonne \\ 31062 TOULOUSE Cedex 9 \\ FRANCE}
\thanks{}

\begin{abstract}
The transfer of the generating operations of an algebra to a homotopy equivalent chain complex produces higher operations. The first goal of this paper is to describe precisely the higher structure obtained when the unary operations commute with the contracting homotopy. To solve this problem, we develop the Koszul duality theory of operads in the category of modules over a cocommutative Hopf algebra. This gives rise to a simpler category of homotopy algebras and infinity-morphisms, which allows us to get a new description of the homotopy category of algebras over such operads. The main example of this theory is given by Batalin--Vilkovisky algebras. 
\end{abstract}

\maketitle

\tableofcontents

\section{Introduction}																

\emph{Homotopy transfer.} In homotopical algebra, one problem is to know how algebraic structures behave under homotopy equivalences. 
For instance, the product of a differential graded associative algebra induces, on a homotopy equivalent chain complex, a product which is not associative in general. However, Kadeishvili proved in \cite{Kadeshvili82} that, in this case, the homotopy equivalent chain complex carries higher operations, in addition to the transferred product. These operations endow it with a homotopy associative algebra structure, also known as $A_{\infty}$-algebra, defined by Stasheff in \cite{Stasheff63}. This is one of the first examples of a \emph{Homotopy Transfer Theorem} (HTT). A HTT was also proved for Lie algebras, for commutative algebras, and more generally, for other types of algebras using the theory of operads, see \cite[Section 10.3]{LodayVallette09}. \\

\emph{Batalin--Vilkovisky algebras.} The theory we develop in this paper is motivated by the example of Batalin--Vilkovisky algebras. These algebras play an important role in geometry, topology and mathematical physics, see for instance \cite{BatVil81,BarKon98,Manin99,ChasSullivan99,Getzler94,LianZuckerman93}. In brief, a BV-algebra structure is given by a commutative product, a Lie bracket and a unary square-zero BV-operator, satisfying some compatibility relations. A HTT for BV-algebras was given in \cite{GalTonVal09}, in terms of a new notion of homotopy BV-algebra. This structure is very rich, and at the same time very intricate. 

However, when the BV-operator commutes with the underlying contracting homotopy, it induces a unary operator which still squares to zero. Thus, there are no higher operations arising from the BV-operator when we apply the HTT of \cite{GalTonVal09}. It is natural to ask what is the precise structure to which it reduces. In this paper, we solve this question as follows: we insert the square-zero operator in the underlying category of chain complexes and, then, we work with operads in this category. This implies that we have to work with operads in the category of modules over a \emph{cocommutative Hopf algebra}. For instance, to encode the category of BV-algebras, we use the operad of Gerstenhaber algebras, enriched with an action of the dual numbers Hopf algebra.\\

\emph{Relative Koszul duality.} In \cite{SalvatoreWahl03}, Salvatore and Wahl began the study of operads with an action of a Hopf algebra. Here, we go even further and we extend the classical Koszul duality theory to this framework in \S\ref{HKoszul}. At each step of the operadic theory, we show that all the objects can be enriched with a compatible action of a Hopf algebra and we prove that the results still hold in this context. In particular, we extend the bar--cobar constructions in \S\ref{barcobar}. We define a notion of homotopy algebras and their infinity-morphisms in \S\ref{Hinfalg}, in order to prove a new Homotopy Transfer Theorem in \S\ref{HTT}. We also extend to this context most of the results, that can be found in \cite{LodayVallette09}, about the homotopy theory of algebras over an operad in \S\ref{Hhoalg}. 

\begin{thmA}\label{ThmA} Let $\bialg$ be a cocommutative Hopf algebra and let $\opd$ be a Koszul operad in the category of $\bialg$-modules. We consider two homotopy equivalent chain complexes in the category of $\bialg$-modules, such that one of them is endowed with a compatible $\opd$-algebra structure. Then the second chain complex is endowed with a $\opd_{\infty}$-algebra structure compatible with the $\bialg$-action.
\end{thmA}

The main example of two homotopy equivalent chain complexes is given by a differential graded algebra and its underlying homology. In this case, being homotopy equivalent in the category of $\bialg$-modules means just that the homotopy and the $\bialg$-action has to commute.

The homotopy category is obtained by localizing the initial category with respect to the quasi-isomorphisms. The $\infty$-morphisms provide a different description of the homotopy category of algebras over an operad. It is easier to deal with this description since the $\infty$-morphisms admit a homotopy inverse, see \cite{Vallette13}. Encoding a category of algebras with an operad in $\bialg$-modules, we get a new description of the associated homotopy category, using simpler $\infty$-morphisms, as follows.

\begin{thmB} Let $\bialg$ be a cocommutative Hopf algebra and let $\opd$ be a Koszul operad in the category of $\bialg$-modules. The homotopy category of dg $\opd$-algebras in $\bialg$-modules is equivalent to the homotopy category of $\opd$-algebras in $\bialg$-modules, together with their $\infty$-$\bialg$-morphisms.    
\end{thmB}
In the case of BV-algebras, we define a simpler category of homotopy BV-algebras: a \emph{strict homotopy BV-algebra} is a ho\-mo\-to\-py Gerstenhaber algebra together with a compatible square-zero unary operator action. Theorem A implies that, when the BV-operator and the contracting homotopy commute, the transferred homotopy BV-algebra structure reduces exactly to a strict homotopy BV-algebra. Theorem B gives a new description of the homotopy category of BV-algebras. We also get a new way to prove the existence of a zigzag of quasi-isomorphisms of BV-algebras, which could help to study the mirror symmetry conjecture by Kontsevich \cite{Kontsevich95}, see \cite{CaoZhou01,CaoZhou03}.\\

\emph{Related literature.} The idea to put unary operations in the underlying category is already present in \cite{GinzburgKapranov94}, but not exactly in the same way as we do. We both put the algebra of \emph{all} the arity one operations in the underlying category, removing them from the operad. In \cite{GinzburgKapranov94}, they keep track of their action using the tensor product over this algebra. Instead, we also remove, from the operad, all the operations made up of the unary ones, and we keep track of the action of these operations via the structure of module over their algebra. Thus, we can still work with the tensor product over the ground field. So, homological results hold in our context, while, in \cite{GinzburgKapranov94}, one needs the algebra of unary operations to be semi-simple. Nevertheless, we need this algebra to form a Hopf algebra. Moreover, since Ginzburg and Kapranov remove \emph{only} the unary operations from the operad, their homotopy algebra structure is much bigger than ours.

There are so far two ways of proving a HTT. On the one hand, one can prove it with model category arguments as in \cite{BergerMoerdijk03}. This relies on a compatibility between the monoidal and the model structures of the underlying category. But, to the best of our knowledge, there is not yet a monoidal model category structure on the category of modules over a cocommutative Hopf algebra together with the tensor product over the ground field, such that the weak equivalences are the quasi-isomorphisms. On the other hand, one can prove a HTT with explicit formulae using the Koszul duality theory for operads, which enables one to prove formality results at the level of algebras for instance. So, we chose this latter method to prove a HTT for operads over Hopf algebras.\\

\emph{Layout}. In the first section, we describe the notion of operads in the category of modules over a Hopf algebra and the associated objects. We give some examples of algebras over such operads. In Section \ref{dualitéKoszul}, we extend the classical Koszul duality to this framework. In Section \ref{Homotopy}, we prove a HTT for algebras over Koszul operads in this context, and we apply it to the example of Batalin--Vilkovisky algebras. Furthermore, we study the homotopy theory of such algebras. \\

\emph{Setting.} Throughout this paper, we work over a field $\KK$ of characteristic zero and all the $\gsym$-modules we consider are reduced, that is the arity zero space is reduced to zero. 

\section{Operads over Hopf algebras }

\subsection{The monoidal category of symmetric modules over a cocommutative Hopf algebra}

We construct the monoidal category of (differential graded) $\gsym$-modules over a cocommutative Hopf algebra, following the example of $\gsym$-modules given in \cite[Section 5.1]{LodayVallette09}. We refer to the book \cite{Sweedler81} for more details on Hopf algebras.

\subsubsection{Recollection on the symmetric monoidal closed category of modules over a Hopf algebra}\index{algèbre de Hopf!module} 

Let $(\bialg,\prdalg,\coprdalg,\unalg,\counalg)$ be a cocommutative bialgebra, where $(\prdalg,\unalg)$ and $(\coprdalg,\counalg)$ are respectively the unital associative algebra structure and the counital cocommutative coalgebra structure. 

\noindent We will denote $1_{\bialg}:=\unalg(1_{\KK})$. Then, we have $\counalg(1_{\bialg})=1_{\KK}$ and $\bialg \cong \textrm{Ker}(\counalg:\bialg \rightarrow \KK) \oplus \KK1_{\bialg}$.\\
We define the iterated coproduct $\coprdalg^n:\bialg \rightarrow \bialg^{\otimes n+1}$ by $\coprdalg^0=\Id$, $\coprdalg^1=\coprdalg$ and $$\coprdalg^n:=(\coprdalg \otimes \Id \otimes \ldots \otimes \Id) \circ \coprdalg^{n-1}.$$
The category $(\HMod, \otimes, \KK)$ of left-$\bialg$-modules is a symmetric monoidal category where the structure of $\bialg$-module on 
\begin{itemize}
 \item[$\bullet$] $\KK$ is defined by the map $\xymatrix{\bialg \otimes \KK \ar[r]^{\counalg \otimes \KK} & \KK \otimes \KK \ar[r] & \KK}$ where the last map is the scalar multiplication,
 \item[$\bullet$] the tensor product of two $\bialg$-modules $M$ and $N$ is defined by
$$ h.(m \otimes n)= \Sw{h}{0}.m \otimes \Sw{h}{1}.n , \forall h \in \bialg, \forall m \in M, \forall n \in N \ ,$$
where $\coprdalg(h)=\Sw{h}{0} \otimes \Sw{h}{1}$.  The Hopf compatibility relation ensures that it is a left $\bialg$-action.
\end{itemize}
We assume moreover that $\bialg$ is a \emph{Hopf algebra}, \index{algèbre de Hopf} that is $\bialg$ is equipped with a linear map $\anti:\bialg \rightarrow \bialg$, called \emph{antipode}, \index{algèbre de Hopf!antipode} which is an inverse of the identity for the convolution product $ \star$: $$ \anti \star \Id = \unalg \counalg = \Id \star \anti \ . $$

Then, for any $\bialg$-modules $M$ and $N$, the vector space $\Hom_{\KK}(M,N)$ carries an $\bialg$-module structure, given by
$$\begin{array}{ccl}
\bialg \otimes \Hom_{\KK}(M,N) & \rightarrow & \Hom_{\KK}(M,N) \\ h \otimes f & \mapsto & h.f : m \mapsto \Sw{h}{0}.f(\anti(\Sw{h}{1}).m) \ .
\end{array}$$

\begin{lem}
The category of $\bialg$-modules is closed.
\end{lem}

\begin{proo}
For any triple $(A,B,C)$ of $\bialg$-modules, the map 
$$ \begin{array}{ccl}
\Hom(A\otimes B,C) & \rightarrow & \Hom(A, \Hom(B,C)) \\
f & \mapsto & (a \mapsto f(a\otimes -))
\end{array}$$
is well defined. It is easy to check that it is a natural isomorphism.
\end{proo}

The previous results extend to the category of differential graded modules.

\begin{pro}\label{symmonclosedcat}
The category $(\dgHMod,\otimes,I)$ of differential graded left-$\bialg$-modules is a symmetric monoidal closed category.
\end{pro}

\begin{Ex}
Consider the algebra $\algnbduaux:= \KK[\unopalg]/(\unopalg^{2})$ of dual numbers\index{algèbre de Hopf!algèbre des nombres duaux}, where $\unopalg$ is of degree $+1$. It is a cocommutative Hopf algebra where the coproduct $\coprdalg:\algnbduaux \rightarrow \algnbduaux \otimes \algnbduaux $ and the antipode $\anti: \algnbduaux \rightarrow \algnbduaux $ are respectively given by 
$$ \coprdalg(\unopalg):=1\otimes\unopalg + \unopalg \otimes 1 $$
and by
$$\anti(\unopalg):= - \unopalg \ .$$
A $\algnbduaux$-module is simply a graded vector space endowed with a map of degree $+1$ that squares to zero, \emph{i.e.} a cochain complex.
\end{Ex}

\subsubsection{The monoidal category of \SH-modules}\label{SHMod} The notion of $\gsym$-modules can be defined in any symmetric monoidal category. We make explicit these objects when the underlying category is equal to the category of $\bialg$-modules.

\begin{defi}
An \emph{\SH-module} \index{algèbre de Hopf!s-module@$\gsym$-module} $M$ is an $\gsym$-module $M=\{ M(n) \}_{n \in \NN}$, such that each $M(n)$ has a left-$\bialg$-module structure which commutes with the $\gsym_n$-module structure.
For $\mu \in M(n)$, the integer $n$ is called the \emph{arity} of $\mu$ and $\mu$ is called an \emph{n-ary operation}.\\
%
A \emph{morphism of \SH-modules} $f:M\rightarrow N$ is a morphism of $\gsym$-modules such that each $f_n:M(n) \rightarrow N(n)$ is $\bialg$-equivariant. \\
We denote the associated category by $\SHMod$. 
\end{defi}

In particular, each $M(n)$ is a $(\bialg,\gsym_{n})$-bimodule and each $f_{n}:M(n) \rightarrow N(n)$ is a morphism of $(\bialg,\gsym_{n})$-bimodules.

\begin{defi}
For any \SH-modules $M$ and $N$, their \emph{tensor product} is the \SH-module $M\otimes N$ defined by 
$$ (M \otimes  N)(n):= \displaystyle{\bigoplus_{i+j=n}} \textrm{Ind}_{\gsym_i \times \gsym_j}^{\gsym_n} M(i) \otimes N(j) \ ,$$
where the action of $\bialg$ is induced by the $\bialg$-module structure on the tensor product of two $\bialg$-modules. Since $\bialg$ is cocommutative, the tensor product of $\bialg$-modules is symmetric.\\
\end{defi}

\begin{pro}\label{Schurtensprd}
The tensor product of \SH-modules is associative with unit the \SH-module $(\KK, 0, 0, \ldots)$. 
\end{pro}

\begin{proo}
By \cite[Proposition 5.1.5]{LodayVallette09}, we have an isomorphism of associativity of the underlying $\gsym$-modules. It is a morphism of $\bialg$-modules by coassociativity of the coproduct. 
\end{proo}

\begin{defi}
For any \SH-modules $M$ and $N$, their \emph{composite product} is the \SH-module $M\circ N$ defined by 
$$ M \circ N:= \displaystyle{\bigoplus_{k \geq 0}} M(k) \otimes_{\gsym_{k}} N^{\otimes k} \ .$$
Here $N^{\otimes k} $ stands for the tensor product of $k$ copies of the \SH-module $N$.\\
For any couple of morphisms of \SH-modules $f:M\rightarrow N$ and $g:M'\rightarrow N'$, their \emph{composite product} is the morphism of \SH-modules $f \circ g: M \circ N \rightarrow M' \circ N'$ given explicitly by the formula:
$$f \circ g (\mu; \nu_1, \ldots, \nu_k):= (f(\mu); g(\nu_1), \ldots, g(\nu_k)).$$
\end{defi}
%
\noindent In arity $n$, we have
$$ (M \circ N)(n):=\displaystyle{\bigoplus_{k \geq 0 }} \ M(k) \otimes_{\gsym_k} \left( \displaystyle{\bigoplus_{i_1+\cdots + i_k=n}}  \textrm{Ind}_{\gsym_{i_1}\times \cdots \times \gsym_{i_k}}^{\gsym_n} N(i_1) \otimes \ldots \otimes N(i_k) \right) \ ,  $$
where the action of $\bialg$ is induced by the $\bialg$-module structure on tensor products of $\bialg$-modules.

\begin{pro}\label{monoidalSHMod}
The category of \SH-modules $(\SHMod, \circ, I)$ is a monoidal category, where the \SH-module $I=(0,\KK, 0, \ldots)$ concentrated in arity $1$ is called \emph{identity} \SH-module.
\end{pro}

\begin{proo}
By \cite[Proposition 5.1.14]{LodayVallette09}, the category of $\gsym$-modules $(\gsym$-$\textsf{Mod},\circ,I)$ is monoidal. Since the composite product of \SH-modules is defined as the composite product of the underlying $\gsym$-modules endowed with the $\bialg$-action induced by the coproduct, the coassociativity and the counitarity of the coproduct of $\bialg$ imply that the isomorphisms of associativity and unit of the underlying $\gsym$-modules are compatible with the $\bialg$-module structure.
\end{proo}

As for $\gsym$-modules, the composite product of two \SH-modules is not linear on the right-hand side. However, the infinitesimal composite product, defined in Section 6.1.1 of~\cite{LodayVallette09}, extends to the category $\SHMod$ to produce a product $\circ_{(1)}$ which is linear on the right-hand side. Recall that the infinitesimal product is defined by
$$ \infprd{M}{N}:= M \circ (I;N) \ ,$$
where, for any $\gsym$-modules $P$, $P_{1}$ and $P_{2}$, $P \circ (P_{1};P_{2})$, is the following $\gsym$-module 
\small$$(P\circ (P_{1};P_{2}))(n):= \displaystyle{\bigoplus_{k = 1 }^{n}} \ P(k) \otimes_{\gsym_k} \left(  \displaystyle{\bigoplus_{i_{1}+\cdots + i_{k}=n} \bigoplus_{j=1}^{k}} \ \textrm{Ind}_{\gsym_{i_{1}}\times \cdots \times \gsym_{i_{k}}}^{\gsym_n}  P_{1}(i_{1})\otimes \ldots \otimes \underbrace{P_{2}(i_{j})}_{\textrm{$j^{th}$ position}} \otimes \ldots \otimes P_{1}(i_{k}) \right) \ .$$
\normalsize
The infinitesimal objects $\infprd{f}{g}$, $f \circ ' g$, $\comp_{(1)}$ and $\decomp_{(1)}$, defined in Section 6.1 of~\cite{LodayVallette09}, also extend to the category $\SHMod$.

\subsubsection{Differential graded framework}

\begin{defi}
A \emph{differential graded \SH-module}, \index{algèbre de Hopf!s-module différentiel gradué@$\gsym$-module différentiel gradué} or dg \SH-module for short, is a dg $\gsym$-module $(M,d)$ such that each $M(n)$ is an $\bialg$-$\gsym_{n}$-bimodule and the differential $d$ is compatible with the $\bialg$-action.\\
A \emph{morphism of differential graded \SH-modules} $f:(M,d_{M})\rightarrow(N,d_{N})$ is a morphism of the underlying dg $\gsym$-modules which is $\bialg$-equivariant. 
We denote by $\dgSHMod$ the category of dg \SH-modules with their morphisms.
\end{defi}

The objects described in the previous section extend to the differential graded framework. However, they now involve signs in their definition. For more details, see \cite[Section 6.2]{LodayVallette09}.\\
We define the \emph{suspension} and the \emph{desuspension} of a graded \SH-module as the suspension and the desuspension of its underlying graded $\gsym$-module respectively, where the action of $\bialg$ is given by the counit of $\bialg$.

\begin{pro}
The category of $(\dgSHMod, \circ , I )$ is a monoidal category.
\end{pro}

\begin{proo}
This follows from the combination of Proposition \ref{monoidalSHMod} and \cite[Proposition 6.2.4]{LodayVallette09}.
\end{proo}

\subsection{Operad over a cocommutative Hopf algebra}\label{Hopd}

In \cite[Section 5.2]{LodayVallette09}, the authors give a monoidal definition of an algebraic operad, that is an operad in the category of vector spaces. More generally, one can define the notion of operad in any monoidal category. Here, we point out the extra structure one gets when one considers operads in the category of \SH-modules instead of just $\gsym$-modules.

\subsubsection{Monoidal definition}

\begin{defi}
An \emph{$\bialg$-operad} \index{algèbre de Hopf!opérade} is a monoid $(\opd, \comp, \unit)$  in the monoidal category of \SH-modules. 
A \emph{morphism $\alpha: \opd \rightarrow \mathcal{Q}$ of $\bialg$-operads} is a morphism of monoids in the category of \SH-modules, that is a morphism of \SH-modules which is compatible with the monoidal structures.
We denote the category of $\bialg$-operads by $\HOpd$.
\end{defi}

\begin{Rq}
Actually, an $\bialg$-operad is nothing but an operad $(\opd,\comp,\unit)$, as defined in \cite[Section 5.2.1]{LodayVallette09}, such that each $\opd(n)$ has an $\bialg$-action which makes $\opd$ into an \SH-module and such that the maps $\comp$ and $\unit$ commute with the $\bialg$-action. Furthermore, a morphism of $\bialg$-operads is a morphism of the underlying operads which commutes with the action of $\bialg$.
\end{Rq}



When $A$ is an $\bialg$-module, the following $\gsym$-module
 $$\End_A(n):=\Hom_{\KK}(A^{\otimes n}, A)$$
is endowed with an $\bialg$-module structure, according to the previous section, which makes it into an \SH-module. This structure is explicitly given by 
$$ \begin{array}{rcl}
    \bialg \otimes \End_A(n) & \rightarrow & \End_A(n)\\
      h \otimes f & \mapsto & h.f:= a_1 \otimes \ldots \otimes a_n  \mapsto  \Sw{h}{0}.f(\anti(\Sw{h}{1}).a_1, \ldots , \anti(\Sw{h}{n}).a_n)
   \end{array} \ .$$

\begin{lem}\label{endopd}
The triple $(\End_{A}$, $\comp_{\End_{A}}$, $\unit_{\End_{A}}: 1_{\KK} \mapsto \Id_{A})$, where $\comp_{\End_{A}}$ is the usual composition map given by $\comp_{\End_{A}}(f;g_{1},\ldots,g_{k}):=f \circ (g_{1}\otimes \ldots \otimes g_{k})$, is an $\bialg$-operad. 
\end{lem}

\begin{proo}
It is immediate to check that $\End_{A}$ is an operad. 
So, we have to prove that $\comp_{End_{A}}$ and $\unit_{\End_{A}}$ are $\bialg$-equivariant maps. For $(f;g_{1},\ldots,g_{k})\in \End_{A}(n)$, $a_{1}\otimes \ldots \otimes a_{n} \in A^{\otimes n}$ and $h \in \bialg$, we have \\
$\comp_{\End_{A}}(h.(f;g_{1},\ldots,g_{k}))(a_{1}\otimes \ldots \otimes a_{n})$
\begin{eqnarray*}
 & = & \Sw{h}{0}. f \left( \displaystyle{\bigotimes_{j=1}^{k}}[\anti(\Sw{h}{n+j)(0})\Sw{h}{n+j)(1}].g_{j}\left(\displaystyle{\bigotimes_{l=1}^{i_{j}}}\anti(\Sw{h}{i_{1}+\cdots+i_{j-1}+l}).a_{i_{1}+\cdots+i_{j-1}+l}\right)\right)\\
 & = & \Sw{h}{0}. f \left( \displaystyle{\bigotimes_{j=1}^{k}}\unalg\counalg(\Sw{h}{n+j}).g_{j}\left(\displaystyle{\bigotimes_{l=1}^{i_{j}}}\anti(\Sw{h}{i_{1}+\cdots+i_{j-1}+l}).a_{i_{1}+\cdots+i_{j-1}+l}\right)\right)\\
 & = & \Sw{h}{0}.(f\circ(g_{1}\otimes \ldots \otimes g_{k}) (\anti(\Sw{h}{1}).a_{1}\otimes \ldots \otimes \anti(\Sw{h}{n}).a_{n})\\
 & = &( h. \comp_{\End_{A}}(f;g_{1},\ldots,g_{k}))(a_{1}\otimes \ldots \otimes a_{n}) \ . 
\end{eqnarray*}
The $\bialg$-equivariance of $\unit_{\End_{A}}$ follows from the $\bialg$-module structure of $\KK$ and from the fact that the antipode is the inverse of the identity for the convolution product.
\end{proo}

Recall that the \emph{Hadamard product $M\hadprd N$} of two $\gsym$-modules $M$ and $N$ is defined to be the following $\gsym$-module 
$$M \hadprd N (n):= M(n) \otimes N(n) \ .$$

\begin{pro}
The Hadamard product of the underlying $\gsym$-modules of two $\bialg$-operads, endowed with the $\bialg$-module structure on the tensor product of two \SH-modules, is an $\bialg$-operad.
\end{pro}

\begin{proo}
Let $\opd$ and $\mathcal{Q}$ be two $\bialg$-operads. By \cite[Section 5.3.3]{LodayVallette09}, we have to prove that the composition map given there is $\bialg$-equivariant. This is the case since it is the composite of $\bialg$-equivariant maps. 
\end{proo}

\begin{defi}
An \emph{ideal} $\Ical$ of an $\bialg$-operad $\opd$ is a sub-\SH-module of $\opd$ such that the composition $\comp_{\opd}(\mu;\nu_{1},\ldots ,\nu_{k})$ is in $\Ical$ as soon as one of the $ \mu, \nu_{1},\ldots ,\nu_{k}$ is in $\Ical$.
The \emph{quotient of an $\bialg$-operad $\opd$ by the ideal $\Ical$} is the $\bialg$-operad $\opd / \Ical$ given by
$$ (\opd / \Ical )(n):= \opd(n) / \Ical(n) \ , \ \forall n \in \NN$$
where the composition map $\comp_{\opd / \Ical}$ is induced by $\comp_{\opd}$. 
\end{defi}

\subsubsection{The free $\bialg$-operad}\label{freeHOp}

\begin{defi}
A \emph{free $\bialg$-operad} \index{algèbre de Hopf!opérade libre} over an \SH-module $M$ is an $\bialg$-operad $\mathcal{F}M$ equipped with a morphism $\unit_{M}: M \rightarrow \mathcal{F}M$ of \SH-modules, which satisfies the following universal property.\\
\\
\hspace*{1cm} Any morphism $f:M \rightarrow \opd$ of \SH-modules, where $\opd$ is an $\bialg$-operad, extends uniquely into a morphism $\tilde{f}:\mathcal{F}M\rightarrow \opd$ of $\bialg$-operads
$$\xymatrix{M \ar[r]^{\unit_{M}} \ar[dr]_{f} & \mathcal{F}M \ar[d]^{\tilde{f}}  & \\ & \opd &\hspace*{-1cm} .} $$
\end{defi}

In other words, the functor $\mathcal{F}: \SHMod \rightarrow \HOpd$ is left adjoint to the forgetful functor $\sqcup$ from $\HOpd$ to $\SHMod$.\\
\\
Recall that the free operad $\freeop M$ over an $\gsym$-module $M$ admits the following realization 
$$ (\freeop M)(n) \cong \displaystyle{\bigoplus_{\tau \in T'(n)}} \tau(M) \ ,$$ 
where $T'(n)$ denotes a set of representatives of isomorphisms classes of $n$-trees.  For any tree $\tau$, the treewise tensor module $\tau(M)$, defined in Section 2.6 of~\cite{Hoffbeck10}, is given by
$$  \tau(M) \cong \displaystyle{\bigotimes_{v \in V(\tau)}} M(|I_{v}|) \ ,$$
where $V(\tau)$ denotes the set of vertices of $\tau$ and where $I_{v}$ denotes the set of entries of the vertex $v$. Then, if $M$ is an \SH-module, we have an $\bialg$-module structure on the treewise tensor module $\tau(M)$, for any tree $\tau$ in $T'(n)$. This action is given by the following map:
$$\begin{array}{ccc}
\bialg \otimes \tau(M) & \rightarrow & \tau(M) \\
h \otimes m_{i_{1}}\ldots m_{i_{|V(\tau)|}} & \mapsto & \Sw{h}{0}.m_{i_{1}} \ldots \Sw{h}{|V(\tau)|-1}.m_{i_{|V(\tau)|}} \ ,
\end{array}$$
where $\coprdalg^{|V(\tau)|-1}(h):=\Sw{h}{0}\otimes\ldots\otimes\Sw{h}{|V(\tau)|-1}$. It amounts to act on the labeling of each vertex of $\tau$ using the coproduct of $\bialg$.

\begin{pro}\label{freeHopd}
Let $M$ be an \SH-module and $\underlying{M}$ be its underlying $\gsym$-module. The free operad $\freeop\underlying{M}$ endowed with the $\bialg$-action described above is the free $\bialg$-operad over $M$.
\end{pro}

\begin{proo}
The composition map $\comp_{\freeop \underlying{M}}$, which corresponds to the grafting of trees, is an $\bialg$-equivariant map. Moreover, if a morphism $f:M \rightarrow \opd$ of $\gsym$-modules is $\bialg$-equivariant then so is $\widetilde{f}:\freeop \underlying{M} \rightarrow \opd$.
\end{proo}

\begin{pro}\label{quotientopd}
Let $M$ be an $\gsym$-module and $R$ be a sub-$\gsym$-module of $\freeop M$. Then, the operad $\freeop M /(R)$, where $(R)$ is the ideal of $\freeop M$ generated by $R$, is an $\bialg$-operad if and only if
\begin{itemize}
\item[$\bullet$] $M$ is an \SH-module
\item[$\bullet$] $R$ is a sub-\SH-module of $\freeop M$.
\end{itemize}
In this case, the $\bialg$-module structure is induced by the action on each vertex using the coproduct of $\bialg$.
\end{pro}

\begin{proo}
The proof of Proposition 4.2 of ~\cite{SalvatoreWahl03} extends to our linear context.
\end{proo}

\subsubsection{Algebra over an $\bialg$-operad}\index{algèbre de Hopf!algèbre sur une opérade} 

\begin{defi}
Let $\opd$ be an $\bialg$-operad. An \emph{algebra over $\opd$}, or for short an $\bialg$-$\opd$-algebra, is an \mbox{$\bialg$-module} $\alg$ equipped with an $\bialg$-equivariant map $\comp_\alg: \opd(\alg)\rightarrow \alg$, which is compatible with the monoidal structure of $\opd$.
A \emph{morphism $f:\alg \rightarrow \alg'$ of $\bialg$-$\opd$-algebras} is a morphism of $\bialg$-modules which commutes with $\comp_{\opd}$.
We denote by $\bialg$-$\opd$-$\mathsf{Alg}$ the category of $\bialg$-$\opd$-algebras.
\end{defi}


If $\alg$ is an $\bialg$-$\opd$-algebra, then the map $\comp_\alg$ is $$\comp_\alg: \opd(A)=\displaystyle{\bigoplus_{n \geq 0}} \ \opd(n) \otimes_{\gsym_n} \alg^{\otimes n} \rightarrow \alg.$$ For $\mu \in \opd(n)$ and $a_1 \otimes \ldots \otimes a_n \in \alg^{\otimes n}$, we denote $\comp_\alg(\mu; a_1, \ldots, a_n)$ simply by $\mu(a_1, \ldots, a_n)$.


\begin{pro}\label{Hopdalg}
Let $\opd$ be an $\bialg$-operad. An $\bialg$-$\opd$-algebra structure on an $\bialg$-module $A$ is equivalent to a morphism of $\bialg$-operads $\opd \rightarrow \End_A$.
\end{pro}

\begin{proo}
By Proposition 5.2.10 of \cite{LodayVallette09}, it remains to prove that $ \comp_{A}:\opd(n)\otimes_{\gsym_{n}} A^{\otimes n} \rightarrow A $ is $\bialg$-equivariant if and only if $\alpha: \opd(n) \rightarrow \End_{A}(n)$ is $\bialg$-equivariant. This follows from the closed structure of the category $\HMod$.
\end{proo}

\begin{pro}[\cite{SalvatoreWahl03}]\label{semidirectprd}
Let $\opd$ be a graded $\bialg$-operad. There exists a graded operad, the semi-direct product $\opd \rtimes \bialg$, such that the category of $\bialg$-$\opd$-algebras, that is $\bialg$-modules with an action of the $\bialg$-operad $\opd$, is isomorphic to the category of $\opd \rtimes \bialg$-algebras, that is modules with an action of the operad $\opd \rtimes \bialg$. 

The operad $\opd\rtimes\bialg$ given by the $\gsym$-module defined by
$$\opd \rtimes \bialg (n):= \opd(n) \otimes \bialg^{\otimes n} \ , $$
where $\gsym_{n}$ acts diagonally acting on $\opd(n)$ and permuting the elements of $\bialg^{\otimes n}$, together with the composition map defined by
$$ \comp_{\opd\rtimes\bialg}((\mu\otimes h);(\nu_{1}\otimes g_{1}),\ldots,(\nu_{k}\otimes g_{k})):=\comp_{\opd}(\mu;h_{(0)}^{1}\cdot\nu_{1},\ldots,h_{(0)}^{k}\cdot\nu_{k})\otimes h_{(1)}^{1}\cdot g_{1}\otimes\cdots\otimes h_{(1)}^{k}\cdot g_{k} \ , $$
where $h=h^{1}\otimes\cdots\otimes h^{k}$, $g_{i}=g_{i}^{1}\otimes\cdots\otimes g_{i}^{n_{i}}$ and where $h^{i}$ acts on $g_{i}$ via the coproduct of $\bialg$. The unit of $\opd\rtimes\bialg$ is given by $\id\otimes 1_{\bialg}$.
\end{pro}


%
%

\subsubsection{Cooperads over a cocommutative Hopf algebra}\index{algèbre de Hopf!coopérade} 

\begin{defi}
An \emph{$\bialg$-cooperad} is a comonoid $(\coopd, \decomp, \coun)$ in the monoidal category of \SH-modules. A morphism of $\bialg$-cooperads $f:\coopd\rightarrow\mathcal{D}$ is a morphism of the underlying \SH-modules compatible with the comonoidal structure of $\coopd$ and $\mathcal{D}$. We denote by $\HCoopd$ the category of $\bialg$-cooperads.
There is an element $\textrm{id}\in \coopd(1)$ such that $\coun(\textrm{id})=1_{\KK}$ and which is called the \emph{identity cooperation}.\\
An $\bialg$-cooperad $\coopd$ is said to be \emph{coaugmented} if there is a morphism $\upsilon:I \rightarrow \coopd$ of $\bialg$-cooperads such that $\coun \upsilon = Id_{I}$.
\end{defi}

The notion of $\bialg$-cooperad is dual to the one of $\bialg$-operad. From the explicit description of the composite product of two \SH-modules and the isomorphism between invariants and coinvariants, it follows that $\decomp$ is made up of $\bialg$-$\gsym_{n}$-bimodules morphisms
$$\decomp(n):\coopd(n) \rightarrow (\coopd \circ \coopd)(n)=  \displaystyle{\bigoplus_{k = 1 }^{n}} \ \coopd(k) \otimes_{\gsym_k} \left( \displaystyle{\bigoplus_{i_1+\cdots + i_k=n}}  \textrm{Ind}_{\gsym_{i_1}\times \cdots \times \gsym_{i_k}}^{\gsym_n} \coopd(i_1) \otimes \ldots \otimes \coopd(i_k) \right) \ .$$


%

\begin{pro}\label{hadprdcoopd}
The Hadamard product of two $\bialg$-cooperads carries an $\bialg$-cooperad structure.
\end{pro}

\begin{proo}
Let $\coopd$ and $\mathcal{D}$ be two $\bialg$-cooperads. By Section 8.3.3 of \cite{LodayVallette09}, it remains to prove that the decomposition map given there is $\bialg$-equivariant. This is the case since this decomposition map is the composite of $\decomp_{\coopd}\otimes \decomp_{\mathcal{D}}$ with maps of the type $\Id \otimes \ldots \otimes \tau \otimes \ldots \otimes \Id$, which are both $\bialg$-equivariant. 
\end{proo}

\begin{defi}
Let $M$ be an \SH-module such that $M(0)=0$. The \emph{cofree $\bialg$-cooperad on $M$} is the $\bialg$-cooperad, denoted by $\freeop^{c}(M)$, which is cofree in the category of conilpotent $\bialg$-cooperads. It means that $\freeop^{c}(M)$ satisfies the following universal property:\\
\\
For any morphism of \SH-modules $\Phi:\coopd \rightarrow M$, from a conilpotent $\bialg$-cooperad $\coopd$, such that $\Phi(\id)=0$, there exists a unique morphism $\widetilde{\Phi}:\coopd \rightarrow \freeop^{c}(M)$ of $\bialg$-cooperads whose corestriction to $M$ is equal to $\Phi$ 
$$ \xymatrix{\coopd \ar[r]^{\Phi} \ar[dr]_{\widetilde{\Phi}} & M \\ & \freeop^{c}(M) \ar@{->>}[u] \ .}  $$
\end{defi}

%

\subsubsection{From operads to cooperads and vice-versa}\label{opdtocoopd}

In the classical case, the aritywise linear dual $\gsym$-module $\coopd^{*}=\{ \Hom_{\KK}(\coopd(n),\KK)\}_{n \in \NN}$ associated to a cooperad $\coopd$ carries an operad structure. The unit is obtained by dualization of the counit $\coun_{\coopd}$ and the composition map is obtained by dualization of $\decomp_{\coopd}$ composed with the natural map from invariants to coinvariants given in Section 5.1.21 of \cite{LodayVallette09}. Since $\bialg$ is a Hopf algebra, if $\coopd$ is moreover an $\bialg$-cooperad then its aritywise dual has an $\bialg$-operad structure. Equally, if $\opd$ is an $\bialg$-operad, such that each $\opd(n)$ is finite dimensional, then its linear dual has an $\bialg$-cooperad structure.


\subsubsection{Differential graded framework}\label{dgframework}\index{algèbre de Hopf!opérade différentielle graduée} 

As for operads, the notion of $\bialg$-operad extends to the differential graded framework. In particular, a \emph{differential graded $\bialg$-operad} is a monoid $(\opd, d_{\opd}, \comp, \unit)$ in the monoidal category $(\dgSHMod,\allowbreak \circ,\allowbreak I)$, that is $(\opd,d_{\opd})$ is a dg \SH-module and $(\opd, \comp, \unit)$ is an $\bialg$-operad structure on $\opd$, such that $\comp$ and $\unit$ are morphisms of dg \SH-modules. Moreover, the results on dg operads of \cite[Section 6.3]{LodayVallette09} can be extended to dg $\bialg$-operads.
Let us recall that we denote by $\suspopd$ (resp. $\suspopd^{-1}$) the endomorphism $\bialg$-operad associated to the graded \SH-module $\KK s$ (resp. $\KK s^{-1}$), see Lemma \ref{endopd}.

\subsection{Examples}

\subsubsection{Mixed chain complexes}


A \emph{mixed chain complex} \index{complexe de chaînes mixte} is a graded vector space $V=V_{\bullet}$ endowed with two maps: a degree $-1$ map $d$ and a degree $+1$ map $\delta$, which both square to zero and anti-commute. One can see any mixed chain complex either as a dg $\algnbduaux$-algebra, or as a dg $\algnbduaux$-$I$-algebra, where $I$ is the identity $\algnbduaux$-operad.


\subsubsection{Batalin--Vilkovisky algebras}

\begin{oldefi}
A \emph{Gerstenhaber algebra} \index{algèbre!Gerstenhaber} is a differential graded vector space $(A,d_{A})$ endowed with
\begin{itemize}
\item[$\diamond$] a symmetric binary product $\prd{}{}$ of degree $0$,
\item[$\diamond$] a symmetric bracket $\bracket{\, }{ }$ of degree $+1$,
\end{itemize}
such that $d_{A}$ is a derivation with respect to each of them and such that
\begin{itemize}
\item[$\rhd$] the product $\prd{}{}$ is associative
$$ (\prd{(\prd{\textrm{-}}{\textrm{-}})}{\textrm{-}}) = (\prd{\textrm{-}}{(\prd{\textrm{-}}{\textrm{-}})})  \ , $$
\item[$\rhd$] the bracket $\bracket{ \, }{  }$ satisfies the Jacobi identity
$$  \bracket{\bracket{\textrm{-} \, }{\textrm{-}}}{\textrm{-}}+\bracket{\bracket{\textrm{-} \,}{\textrm{-}}}{\textrm{-}}\cdot(123)+\bracket{\bracket{\textrm{-} \, }{\textrm{-}}}{\textrm{-}} \cdot(321)=0 \ , $$
\item[$\rhd$] the product $\prd{}{}$ and the bracket $\bracket{ \ }{ \, }$ satisfy the Leibniz relation
$$ \bracket{\textrm{-} \, }{\prd{ \textrm{-}}{\textrm{-}}} = (\prd{\bracket{\textrm{-} \, }{\textrm{-}}}{\, \textrm{-}})+ (\prd{\textrm{-} \, }{\bracket{\textrm{-} \, }{\textrm{-}}} ) \cdot (12) \ . $$
\end{itemize}
\indent A \emph{Batalin--Vilkovisky algebra}, \index{algèbre!Batalin--Vilkovisky} or \emph{BV-algebra} for short, is a Gerstenhaber algebra $A$ endowed, in addition, with 
\begin{itemize}
\item[$\diamond$] a unary operator $\unop$ of degree $+1$ ,
\end{itemize}
such that $d_{A}$ is a derivation with respect to $\unop$ and such that
\begin{itemize}
\item[$\blacktriangleright$] the operator satisfies $\unop^{2}=0$, 
\item[$\blacktriangleright$] the bracket is the obstruction of $\unop$ being a derivation with respect to the product $\prd{}{}$
$$  \bracket{\textrm{-} \, }{\textrm{-}} = \unop \circ (\prd{\textrm{-}}{\textrm{-}}) - (\prd{\unop(\textrm{-})}{\textrm{-}})-(\prd{\textrm{-}}{\unop(\textrm{-})})  \ , $$
\item[$\blacktriangleright$] the operator $\unop$ is a graded derivation with respect to the bracket $\bracket{}{}$
$$  \unop \circ (\bracket{\textrm{-} \, }{\textrm{-}})+\bracket{\unop(\textrm{-})}{\textrm{-}}+\bracket{\textrm{-} \, }{\unop(\textrm{-})} = 0 \ . $$
\end{itemize}
\end{oldefi}

Let $\Gerst$ be the operad encoding Gerstenhaber algebras. It is defined by generators and relations as follows
$$ \Gerst := \freeop(\gensp)/(\rel) \ ,$$
where $\gensp= E(2):=\prd{}{} \KK_{2} \oplus \bracket{ \, }{} \KK_{2}$, with $\KK_{2}$ the trivial representation of the symmetric group $\gsym_{2}$. The space of relations $\rel$ is the sub-$\gsym$-module of $\freeop(\gensp)$ generated by the relations $\rhd$.
The operad $\BV$, encoding BV-algebras, is then given by
$$ \freeop(\gensp')/(R') \ ,$$
where $\gensp':=\gensp \oplus \unop\KK$ and where $R'$ is the sub-$\gsym$-module of $\freeop(\gensp')$ generated by the relations $\rhd$ and $\blacktriangleright$.\\

\begin{lem}
The action of $\algnbduaux$ on the generators $\prd{}{}$ and $\bracket{\,}{}$ of $\Gerst$ given by
$$ \begin{array}{cccc}
\unopalg  : & \prd{}{} & \mapsto & \bracket{\,}{}\\
& \bracket{\,}{} & \mapsto & 0 
\end{array}$$
induces a $\algnbduaux$-operad structure on $\Gerst$.
\end{lem}

\begin{proo}
Using Proposition \ref{quotientopd}, we just check that $\gensp$ is an $\gsym$-$\algnbduaux$-module and that $\rel$ is a sub-$\gsym$-$\algnbduaux$-module of $\freeop(\gensp)$.
\end{proo}

\begin{Rq}
The $\unopalg$-action on the generators of $\Gerst$ corresponds to the effect in homology of the action of the circle group $S^{1}$ on the little discs operad given in~\cite[Section 4]{Getzler94bis}. This structure induced in homology is made explicit in~\cite[Section 5]{SalvatoreWahl03}. 
\end{Rq}

As a consequence of~\cite[Theorem 2.7, Example 4.4]{Markl96D}, we have $\Gerst \cong \Com \circ \Lie_{1}$ as operads. This isomorphism gives us a way to make explicit the $\algnbduaux$-operad structure on $\Gerst$. We 
denote by $\odot$ the commutative tensor product, that is the quotient of the tensor product under the permutation of terms. In particular, a generic element of $\Com \circ \Lie_{1}$ is of the form $L_{1} \odot \cdots \odot L_{t}$, with $L_{i}\in \Lie_{1}$, for $i=1, \ldots, t$; the elements of $\Com$ being implicit.
\begin{pro}
Under the isomorphism of operads $\Gerst \cong \Com \circ \Lie_{1}$,  the structure of $\algnbduaux$-operad on the right-hand side is given by
$$ \unopalg \cdot  (L_{1} \odot \cdots \odot L_{t}) = \displaystyle{\sum_{1\leq i < j \leq n}} (-1)^{\varepsilon_{i,j}} \bracket{L_{i}}{L_{j}}\odot L_{1} \odot \cdots \odot \hat{L}_{i} \odot \cdots \odot \hat{L}_{j} \odot \cdots \odot L_{t} \ , $$
where $\bracket{L_{i}}{L_{j}}:= \comp_{\Lie_{1}}(\bracket{\,}{}; L_{i},L_{j})$ and where the sign $\varepsilon_{i,j}$, arising from the Koszul sign rule, is given by
$$\varepsilon_{i,j}= (|L_{i}|+|L_{j}|)(|L_{1}|+\cdots+|L_{i-1}|) + |L_{j}|(|L_{i+1}|+\cdots+|L_{j-1}|) \ . $$
\end{pro}

\begin{proo}
By induction on $n\geq 2$.
\end{proo}

\begin{Rq}
The action of $\unopalg$ on an element $ L_{1} \odot \cdots \odot L_{t} $ in $\Com \circ \Lie_{1}$ is exactly the image of $\unopalg^{*}\otimes L_{1} \odot \cdots \odot L_{t}$ under the derivation $^{t} d_{\varphi}$, which is given in~\cite[Proof of Lemma 5]{GalTonVal09}. In particular, $\unopalg$ acts as the Chevalley-Eilenberg boundary map defining the homology of Lie algebras.
\end{Rq}

\begin{pro}\label{BValg}
The category of $\algnbduaux$-$\Gerst$-algebras is isomorphic to the category of $\BV$-algebras. 
\end{pro}

\begin{proo}
It is a straightforward consequence of Proposition \ref{semidirectprd} conjugated with the following isomorphism of operads proved by Salvatore and Wahl in \cite{SalvatoreWahl03}: 
$$ \BV \cong \Gerst \rtimes \algnbduaux \ . $$
 \end{proo}

\begin{Rq} The algebraic structure of a Batalin--Vilkovisky algebra can be encoded into two different ways, using the operad $\BV$ or the $\algnbduaux$-operad $\Gerst$. When using the graded $\algnbduaux$-operad $\Gerst$, the unary operator $\unop$ is provided by the underlying category of mixed chain complexes and its relations with the product and the Lie bracket are encoded in the action of $\unopalg$ on those generating operations.
\end{Rq}

More generally, we can consider algebras over the homology of the framed little $n$-discs operad $f\mathcal{D}_{n}$. Batalin--Vilkovisky algebras correspond to the case $n=2$, by \cite{Getzler94}. In this case, Poisson $n$-algebras, that is algebras over the homology of the little $n$-discs operad $\mathcal{D}_{n}$, play the role of Gerstenhaber algebras. In the same way, it is proved in \cite{SalvatoreWahl03} that
$$f\mathcal{D}_{n} \cong \mathcal{D}_{n} \rtimes SO(n) \ , $$
as topological operads. Hence, taking the homology, we obtain that
$$ H(f\mathcal{D}_{n}) \cong H(\mathcal{D}_{n}) \rtimes H(SO(n)) \ . $$
In other words, a structure of algebra over $H(f\mathcal{D}_{n})$ is equivalent to a Poisson $n$-algebra structure in the category of $H(SO(n))$-modules.

However, when $n>3$, the algebra $H(SO(n))$ has more than one generator. So, we can chose to keep some in the operad and to put only a part of these generators in the underlying category. We still have a Hopf algebra, which acts on the operad thus obtained. 

\begin{Rq} Unlike \cite{GinzburgKapranov94}, we can chose to put only a part of the unary operations in the underlying category.
\end{Rq}

\subsubsection{Differential graded algebras over an operad}

Let $\opd$ be an operad. A differential graded $\opd$-algebra is a chain complex $(A,d_{A})$ endowed with a $\opd$-algebra structure such that $d_{A}$ is a derivation for all the operations of $\opd$. Thus, we can encode this structure in the $\algnbduaux$-operad $\opd$, where the action of $\algnbduaux$ is trivial.

\section{Relative Koszul duality theory}\label{dualitéKoszul}

In this section, we settle the Koszul duality theory for operads over a cocommutative Hopf algebra $\bialg$, following the method of \cite[Chapter 7]{LodayVallette09}. At each step, we prove that the objects can be enriched with an action of $\bialg$, and we show that the results still hold in this context. 

Let $(\opd,\comp,\unit,d_{\opd})$ be a dg $\bialg$-operad and $(\coopd,\decomp,\coun,d_{\coopd})$ be a dg $\bialg$-cooperad.

\subsection{Twisting morphisms}\index{algèbre de Hopf!morphisme tordant} 

We consider the following dg \SH-module 
$$(\Hom(\coopd,\opd):=\left\lbrace \Hom_{\KK}(\coopd(n),\opd(n)) \right\rbrace_{n \geq 0} ,\partial) \ ,$$
where $\gsym_{n}$ acts by conjugation and where $\partial (f) = [d_{\opd}, f] := d_{\opd} \circ f - (-1)^{|f|} f \circ d_{\coopd}$.

\begin{pro}
The dg \SH-module $(\Hom(\coopd,\opd),\partial)$ is a dg $\bialg$-operad, called the \emph{convolution $\bialg$-operad}.\index{algèbre de Hopf!opérade de convolution} 
\end{pro}

\begin{proo}
By \cite[Section 1]{BergerMoerdijk03}, the $\gsym$-module $(\Hom(\coopd,\opd),\partial)$ is an operad with respect to the following composition map :
$$ \comp_{\Hom(\coopd,\opd)}(f;g_{1},\ldots,g_{k}) := \comp_{\opd}\circ (f\otimes g_{1}\otimes\cdots\otimes g_{k})\circ \decomp_{\coopd} \ . $$
Moreover, in \cite[Section 6.4.1]{LodayVallette09}, it is proved that $(\Hom(\coopd,\opd),\comp_{\Hom(\coopd,\opd)},\partial)$ is a dg operad.
The map $\comp_{\Hom(\coopd,\opd)}$ is $\bialg$-equivariant as a composite of $\bialg$-equivariant map. And, since $d_{\opd}$ and $d_{\coopd}$ are $\bialg$-equivariant, then so is the differential $\partial$.
\end{proo}


Recall that the following composite defines a \emph{pre-Lie product} on $\displaystyle{\prod_{n \geq 0}} \Hom(\coopd,\opd)(n)$
$$f \star g:= \vcenter{\xymatrix{\coopd \ar[r]^{\hspace{-0.5cm}\decomp_{(1)}} & \infprd{\coopd}{\coopd} \ar[r]^{\hspace{-0.2cm}\infprd{f}{g}} & \infprd{\opd}{\opd} \ar[r]^{\hspace{0.4cm}\comp_{(1)}} & \opd }} \ . $$


We denote by $\Hom_{\bialg\textrm{-}\gsym_{n}}(\coopd(n),\opd(n))$ the space of $\bialg$-$\gsym_{n}$-bimodules morphisms from $\coopd(n)$ to $\opd(n)$ and we denote the associated product of $\gsym$-$\bialg$-equivariant maps by
$$\Hom_{\gsym\textrm{-}\bialg}(\coopd,\opd):= \displaystyle{\prod_{n \geq 0}} \ \Hom_{\bialg\textrm{-}\gsym_{n}}(\coopd(n),\opd(n)) \ . $$

\begin{pro}
The space $(\Hom_{\gsym\textrm{-}\bialg}(\coopd,\opd),\star,\partial)$ is a dg pre-Lie algebra. The associated Lie bracket induces a dg Lie algebra structure $(\Hom_{\gsym\textrm{-}\bialg}(\coopd,\opd),[\ , \ ],\partial)$.
\end{pro}

\begin{proo}
By Proposition 6.4.5 and Lemma 6.4.6 of \cite{LodayVallette09}, it only remains to prove that if $f$ and $g$ are in $\Hom_{\gsym\textrm{-}\bialg}(\coopd,\opd)$ then so is $f\star g$. This is the case since the pre-Lie product $\star$ is defined to be the composite of $\bialg$-equivariant maps.
\end{proo}

We consider the \emph{Maurer-Cartan equation} in the dg-pre-Lie algebra $(\Hom_{\gsym\textrm{-}\bialg}(\coopd,\opd),\star,\partial)$
$$\partial(\twm)+\twm \star \twm=0 \ .$$

\begin{defi}
A solution $\twm:\coopd \rightarrow \opd$ of degree $-1$ to the Maurer-Cartan equation is called a \emph{twisting $\bialg$-morphism}. We denote by $\Tw_{\bialg}(\coopd,\opd)$ the space of twisting $\bialg$-morphisms from $\coopd$ to $\opd$.\\
When $\coopd$ is a coaugmented dg cooperad and $\opd$ is an augmented dg operad, we require that the composition of a twisting morphism with respectively the coaugmentation map or the augmentation map vanishes.
\end{defi}

\subsection{Bar and cobar constructions}\label{barcobar} 

Are the two functors $\Tw_{\bialg}(\coopd,-)$ and $\Tw_{\bialg}(-,\opd)$ representable?\\

Recall from \cite[Section 6.5]{LodayVallette09}, that the \emph{bar construction} of an augmented dg operad $\opd$ is the dg conilpotent cooperad $$\Barc\opd:=(\freeop^{c}(s\overline{\opd}),d_{\Barc\opd}) \ ,$$ 
and that the \emph{cobar construction} of a coaugmented dg cooperad $\coopd$ is the dg augmented operad $$\Cobc\coopd:=(\freeop(s^{-1}\overline{\coopd}),d_{\Cobc\coopd}) \ . $$

\begin{pro}
Let $\opd$ be an $\bialg$-operad and $\coopd$ be an $\bialg$-cooperad. Then, the $\bialg$-module structure of the free operad over an \SH-module makes $\Barc \opd$ and $\Cobc \coopd$ into a dg $\bialg$-cooperad and a dg $\bialg$-operad respectively.
\end{pro}

\begin{proo}
By Proposition \ref{freeHopd}, the free operad $\freeop(s^{-1}\overline{\coopd})$ is an $\bialg$-operad. The differential $d_{\Cobc \coopd}$ is equal to the sum $d_{1}+d_{2}$, where $d_{1}$ is the unique derivation which extends 
$$ \varphi : \xymatrix@M=3pt@C=40pt{s^{-1}\overline{\coopd} \ar[r]^-{\id\otimes d_{\coopd}} & \freeop(s^{-1}\overline{\coopd})} \ , $$
and where $d_{2}$ is the unique derivation which extends 
\begin{flushleft}
$ \psi : \xymatrix@M=3pt@C=40pt{(\KK s^{-1} \otimes \overline{\coopd}) \ar[r]^-{\decomp_{s}\otimes\decomp_{(1)}} & (\KK s^{-1} \otimes\KK s^{-1})\otimes(\overline{\coopd}\circ_{(1)}\overline{\coopd}) \ar[r]^-{\Id\otimes\tau\otimes\Id} & } $
\end{flushleft}
\begin{flushright}
$(\KK s^{-1} \otimes \overline{\coopd}) \circ_{(1)}(\KK s^{-1} \otimes \overline{\coopd}) \cong \freeop(s^{-1}\overline{\coopd})^{(2)} \rightarrowtail \freeop(s^{-1}\overline{\coopd}) \ . $
\end{flushright}
The Proposition 6.3.15 of \cite{LodayVallette09} extends to graded \SH-modules. Indeed, if $\gensp$ is an \SH-module and $\alpha:\gensp \rightarrow \freeop(\gensp)$ is a morphism of \SH-modules, the unique derivation which extends $\alpha$ is a composite of $\bialg$-equivariant maps. Thus, since $\varphi$ and $\psi$ are morphisms of \SH-modules, the derivations $d_{1}$ and $d_{2}$, and hence $d_{\Cobc \coopd}$, are $\bialg$-equivariant maps. \\
In the same way, we prove that $\Barc \opd$ is a dg $\bialg$-cooperad.
\end{proo}

%

\begin{pro}\label{adjunction}
The bar and cobar constructions form a pair of adjoint functors
$$ \Cobc : \{ \text{ conil. dg $\bialg$-cooperads } \} \rightleftharpoons \{ \text{ aug. dg $\bialg$-operads } \} : \Barc \ , $$
such that the adjunction is given by the set of twisting $\bialg$-morphisms. That is for every augmented dg $\bialg$-operad $\opd$ and every conilpotent dg $\bialg$-cooperad $\coopd$, there exist natural isomorphisms
$$\Hom_{\mathsf{dg} \ \HOpd}(\Cobc\coopd,\opd) \cong \Tw_{\bialg}(\coopd,\opd) \cong \Hom_{\mathsf{dg} \ \HCoopd}(\coopd,\Barc\opd) \ . $$
\end{pro}

\begin{proo}
It is a consequence of \cite[Theorem 6.5.10]{LodayVallette09} and of the definitions of the free $\bialg$-operad and the cofree $\bialg$-cooperad over an \SH-module.
\end{proo}

We denote by $\unit_{\coopd}:\coopd \rightarrow \Cobc\coopd$ and $\pi_{\opd} :\Barc\opd \rightarrow \opd$ the universal twisting $\bialg$-morphisms corresponding respectively to the unit $\upsilon:\coopd \rightarrow \Barc\Cobc\coopd$ and to the counit $\epsilon:\Cobc\Barc\opd \rightarrow \opd$ of the adjunction.

\begin{pro}
Any twisting $\bialg$-morphism $\twm : \coopd \rightarrow \opd$ factorizes uniquely through the universal twisting morphisms $\pi_{\opd}$ and $\unit_{\coopd}$ as follows
$$ \xymatrix{ & \Cobc\coopd \ar@{-->}[dr]^{g_{\twm}} & \\ \coopd \ar[ur]^{\unit_{\coopd}} \ar[rr]^{\twm} \ar@{-->}[dr]_{f_{\twm}} & & \opd \ ,  \\ & \Barc\opd \ar[ur]_{\pi_{\opd}}    }  $$
where $g_{\twm}$ is a morphism of dg $\bialg$-operads and $f_{\twm}$ is a morphism of dg $\bialg$-cooperads. 
\end{pro}

\begin{proo}
It is a consequence of the adjunction given in Proposition \ref{adjunction}.
\end{proo}

When dealing with operads in the monoidal category of dg $\gsym$-modules, a twisting morphism $\twm$ in $\Tw(\coopd,\opd)$ is called \emph{Koszul} when a certain chain complex, called the \emph{Koszul complex}, is acyclic, see \cite[Section 6.6.1]{LodayVallette09}. So, the property for a twisting morphism to be Koszul is a homological property and only depends on the differential structures of $\coopd$ and $\opd$. 

\begin{defi}
A twisting $\bialg$-morphism $\twm \in \Tw_{\bialg}(\coopd,\opd)$ is said to be \emph{Koszul} when, seen as a morphism of dg $\gsym$-modules, it is a Koszul morphism.
\end{defi}


By the same argument, Theorem 6.6.2 of \cite{LodayVallette09} extends to $\bialg$-operads to give the following result.

\begin{thm}\label{koszulqi}
Let $\opd$ be a connected weight graded dg $\bialg$-operad and let $\coopd$ be a connected weight graded dg $\bialg$-cooperad. Let $\twm : \coopd \rightarrow \opd$ be a twisting $\bialg$-morphism. The following assertions are equivalent:
\begin{enumerate}
\item the twisting morphism $\twm$ is Koszul,
\item the morphism of dg $\bialg$-cooperads $f_{\twm}:\coopd \rightarrow \Barc\opd$ is a quasi-isomorphism,
\item the morphism of dg $\bialg$-operads $g_{\twm}:\Cobc\coopd \rightarrow \opd$ is a quasi-isomorphism.
\end{enumerate}
\end{thm}

\begin{thm}
The counit $\upsilon:\Cobc\Barc\opd \rightarrow \opd$ and the unit $\epsilon: \coopd \rightarrow \Barc\Cobc\coopd$ of the adjunction are $\bialg$-quasi-isomorphisms of dg $\bialg$-operads and dg $\bialg$-cooperads respectively.
\end{thm}

The resolution $\Cobc\Barc\opd$ is called the \emph{bar-cobar resolution} of $\opd$.

\subsection{Koszul duality of $\bialg$-operads}\label{HKoszul}


\begin{defi}
A \emph{quadratic data} is a pair $(\gensp,\rel)$ made up of a graded \SH-module $\gensp$ and  a graded sub-\SH-module $\rel$ of $\freeop(\gensp)^{(2)}$ called respectively the \emph{generating operations} and the \emph{relations}.\\
The \emph{quadratic $\bialg$-operad} $\opd(\gensp,\rel)$ associated to a quadratic data $(\gensp,\rel)$ is the quotient $\bialg$-operad of $\freeop(\gensp)$, which satisfies the following universal property: for any $\bialg$-operad $\opd$ of $\freeop(\gensp)$ such that the following composite of morphisms of \SH-modules is trivial
$$ \rel \rightarrowtail \freeop(\gensp) \twoheadrightarrow \opd \ , $$
there exists a unique morphism of $\bialg$-operads which makes the following diagram commutative
$$ \xymatrix@M=6pt{ \rel \; \; \ar@{>->}[r]  & \freeop(\gensp) \ar[r] \ar[rd] & \opd \\  & & \opd(\gensp,\rel) \ar[u] \ .  } $$
Dually, we define the \emph{quadratic $\bialg$-cooperad} $\coopd(\gensp,\rel)$ associated to a quadratic data $(\gensp,\rel)$.
\end{defi}

%
%
\begin{pro}\label{quadopd}
Let $(\gensp,\rel)$ be a quadratic data. The quadratic $\bialg$-operad (resp. $\bialg$-cooperad) associated to $(\gensp,\rel)$ is given by the quadratic operad (resp. cooperad) associated to $(\gensp,\rel)$ in the category of \ $\gsym$-modules, endowed with the $\bialg$-module structure induced by the one on the \SH-module $\freeop(\gensp)$.
\end{pro}

\begin{proo}
We only prove this result for the quadratic $\bialg$-operad. In the category of operads, $\opd(\gensp,\rel)$ is given by the quotient $\freeop(\gensp)/(\rel)$ of the free operad on $\gensp$ by the ideal generated by $\rel$. By Proposition \ref{quotientopd}, $\freeop(\gensp)/(\rel)$ is an $\bialg$-operad and, by assumptions on $\rel$ and on $\opd$,  the morphism of operads $\opd(\gensp,\rel) \rightarrow \opd$ turns out to be an $\bialg$-equivariant map. 
\end{proo}

At this point, we can extend the definition of the Koszul dual (co)operad, given in \cite[Section 7.2]{LodayVallette09}, to the category of \SH-modules. In the classical case,  recall that the \emph{Koszul dual cooperad} of quadratic operad $\opd$ is defined by the following cooperad
$$\opd^{\textrm{!`}}:= \coopd(s\gensp,s^{2}\rel) \ ,$$
and the \emph{Koszul dual operad} of $\opd$ by the following operad
$$\opd^{!}:=(\suspopd^{c}\hadprd \opd^{\textrm{!`}})^{*} \ .$$


\begin{pro}
If $\opd$ is a quadratic $\bialg$-operad, then the Koszul dual cooperad (resp. operad) of $\opd$ is an $\bialg$-cooperad (resp. $\bialg$-operad).
\end{pro}

\begin{proo}
By definition, $\opd^{\textrm{!`}}$ is an $\bialg$-cooperad. Since $\decomp_{\End_{\KK s}}$ is $\bialg$-equivariant, the $\bialg$-operad structure on $\opd^{!}$ follows from Proposition \ref{hadprdcoopd} and Section \ref{opdtocoopd}.
\end{proo}

\begin{pro}
Let $(\gensp,\rel)$ be a quadratic data such that $\gensp$ is a finite dimensional \SH-module. The Koszul dual operad $\opd^{!}$ of the quadratic operad $\opd=\opd(\gensp,\rel)$ admits the following quadratic presentation:
$$ \opd^{!} \cong \opd(s^{-1}\mathcal{S}^{-1}\hadprd \gensp^{*},\rel^{\perp}) \ ,$$
where $\rel^{\perp}$ denotes the sub-\SH-module obtained by proper suspension of the operations indexing the vertices of the trees of the orthogonal module $(s^{2}\rel)^{\perp} \subset \freeop(s^{-1}\gensp^{*})^{(2)}$, which is the image of $(s^{2}\rel)^{*}$ under the isomorphism $ (\freeop(\gensp)^{(2)})^{*}\cong \freeop(s^{-1}\gensp^{*})^{(2)} \ . $
\end{pro}

\begin{proo}
By \cite[Proposition 7.2.4]{LodayVallette09}, this result is true in the category of operads. It is left to prove that the isomorphism is an $\bialg$-equivariant map: the restriction of this isomorphism to the space of generators $s^{-1}\mathcal{S}^{-1}\hadprd \gensp^{*}$ is $\bialg$-equivariant then so is the isomorphism. 
\end{proo}

\begin{pro}
Let $(\gensp,\rel)$ be a quadratic data and $\opd=\opd(\gensp,\rel)$ be its associated quadratic $\bialg$-operad. The natural $\bialg$-cooperad inclusion $i:\opd^{\textrm{!`}} \hookrightarrow \Barc\opd$ induces an isomorphism of graded $\bialg$-operads 
$$ i:\opd^{\textrm{!`}} \xrightarrow{\cong}  H^{0}(\Barc^{*}\opd) \ ,$$
where $H^{\bullet}(\Barc^{*}\opd)$ denotes the cohomology groups of the syzygy degree cochain complex associated to $\Barc \opd$, see \cite[Section 7.3.1]{LodayVallette09}.
\end{pro}

\begin{proo}
By \cite[Proposition 7.3.2]{LodayVallette09}, we have an isomorphism of operads and the inclusion $\opd^{\textrm{!`}} \hookrightarrow \Barc^{0}\opd$ is exactly the kernel of the differential defining $H^{\bullet}(\Barc^{*}\opd)$, that is $H^{0}(\Barc^{*}\opd)$. Since this differential is $\bialg$-equivariant, this isomorphism is compatible with the $\bialg$-module structure.
\end{proo}

\begin{lem}
Let $(\gensp,\rel)$ be a quadratic data. Then, the composite $$\kappa: \coopd(s\gensp, s^{2}\rel) \twoheadrightarrow s\gensp \xrightarrow{s^{-1}} \gensp \hookrightarrow \opd(\gensp,\rel) \ , $$ associated to $(\gensp,\rel)$, is a twisting $\bialg$-morphism.
\end{lem}

\begin{proo}
It is straightforward from Lemma 7.4.2 of \cite{LodayVallette09}, since $\kappa$ is a composite of $\bialg$-equivariant maps.
\end{proo}

\begin{thm}[Koszul criterion]
Let $(\gensp,\rel)$ be a quadratic data. The following propositions are equivalent:
\begin{enumerate}
\item[$\mathrm{(1)}$] the inclusion $i:\opd^{\textrm{!`}} \rightarrowtail \Barc\opd$ is a quasi-isomorphism of dg $\bialg$-cooperads,
\item[$\mathrm{(2)}$] the projection $p:\Cobc\coopd\twoheadrightarrow\coopd^{\textrm{!`}}$ is a quasi-isomorphism of dg $\bialg$-operads.
\end{enumerate}
When these propositions hold, we say that $\opd=\opd(\gensp,\rel)$ is a \emph{Koszul operad}.
\end{thm}

\begin{proo}
This follows from Theorem \ref{koszulqi}.
\end{proo}

\subsection{Example}

%
Let us now focus on the example of Batalin--Vilkovisky algebras.\\

In~\cite{GetzlerJones94}, it was proved that the operad $\Gerst$ is Koszul. Then, by the present definition, so is the $\algnbduaux$-operad $\Gerst$. 
\begin{lem}
The action of $\algnbduaux$ on $\Gerst^{\textrm{!`}}$ is induced by
$$\begin{array}{cccc}
\unopalg\cdot \quad : & \Gerst^{\textrm{!`}} & \twoheadrightarrow & \KK s \; \! \bullet{}{} \oplus \KK s \bracket{\,}{} \\
& s \; \! \bullet{}{} & \mapsto & - s \bracket{\,}{} \\
& s \bracket{\,}{} & \mapsto & 0 \ .
\end{array} $$
\end{lem}

\begin{proo}
By definition, the action of $\algnbduaux$ is characterized by its value on the cogenerators of $\Gerst^{\textrm{!`}}$, which is deduced from the one on $\Gerst$. 
\end{proo}

\begin{pro}\label{explicitDcoopstructure}
Under the isomorphism of cooperads $\Gerst^{\textrm{!`}} \cong \suspopd^{-1}\Com_{1}^{c} \circ \suspopd^{-1}\Lie^{c}$, the structure of $\algnbduaux$-cooperad on the right-hand side is given by
$$ \unopalg \cdot (L_{1} \odot \cdots \odot L_{t})= \displaystyle{\sum_{k=1}^{n}} (-1)^{\varepsilon_{k}} L_{1} \odot \cdots \odot L'_{k}\odot L''_{k}\odot  \ldots \odot L_{t} \ , $$
where $L'_{k}\odot L''_{k}$ is the sumless Sweedler's notation for the image of $L_{k}$ under the binary part of $\decomp_{\suspopd^{-1}\Lie^{c}}$, that is  
$\decomp_{\suspopd^{-1}\Lie^{c}} : \suspopd^{-1}\Lie^{c}\rightarrow \suspopd^{-1}\Lie^{c}(2)\otimes (\suspopd^{-1}\Lie^{c}\otimes \suspopd^{-1}\Lie^{c}) $ 
and 
$$ \varepsilon_{k}= |L_{k}|+\cdots+|L_{t}|+1 \ . $$
\end{pro}

\begin{proo}
To make explicit the action of $\unopalg$, we use the $\algnbduaux$-operad structure on $(\Gerst^{\textrm{!`}})^{*} \cong \suspopd \Com_{-1} \circ \suspopd \Lie$, which is given by the isomorphism of operads $(\Gerst^{\textrm{!`}})^{*} \cong \suspopd^{2} \Gerst \cong \suspopd^{2} (\Com \circ \Lie_{1})$.
We have $$\unopalg \cdot (L_{1} \odot \cdots \odot L_{t})= \unopalg \cdot (L_{1}^{*} \odot \cdots \odot L_{t}^{*})^{*} \ , $$ with $L_{1}^{*} \odot \cdots \odot L_{t}^{*} \in \suspopd \Com_{-1} \circ \suspopd \Lie$. The only elements in $\suspopd \Com_{-1} \circ \suspopd \Lie$ whose image under $\unopalg \cdot (L_{1}^{*} \odot \cdots \odot L_{t}^{*})^{*}$ is non-zero are of the form
$$ L_{1}^{*} \odot \cdots \odot (L'_{k})^{*} \odot (L''_{k}) ^{*}\odot \cdots \odot L_{t}^{*} \ , k \in \{ 1 , \ldots , t \} \ ,$$
where $(L'_{k})^{*}, (L''_{k}) ^{*}$ are elements of $ \suspopd \Lie$ such that $L_{k}^{*}=\comp_{\suspopd \Lie}(s^{-1}[\,;];(L'_{k})^{*}, (L''_{k}) ^{*} )$. This image is equal to 
$(-1)^{\varepsilon_{k}+1}$, where $\varepsilon_{k}=|L_{k}^{*}|+\cdots + |L_{t}^{*}|$.
\end{proo}

\begin{Rq}
Dually to the $\algnbduaux$-operad structure on $\Gerst$, the $\algnbduaux$-cooperad structure on $\Gerst^{\textrm{!`}}$ is equal, up to sign, to the image of an element in $\KK[\unopalg]_{1} \circ \suspopd^{-1}\Com_{1}^{c} \circ \suspopd^{-1}\Lie^{c}$ under the coderivation $d_{\varphi}$, given in~\cite[Lemma 5]{GalTonVal09}. 
\end{Rq}

\section{Homotopy algebras and transfer theorem}\label{Homotopy}

In this section, we extend some definitions and results of \cite[Section 10]{LodayVallette09} to the category of $\bialg$-modules, requiring in addition the compatibility with the $\bialg$-module structure. 

\subsection{The category of homotopy algebras}\label{Hinfalg}

Let $\opd$ be a Koszul $\bialg$-operad. By definition, a \emph{homotopy $\bialg$-$\opd$-algebra}, or an \emph{$\bialg$-$\opd_{\infty}$-algebra}, is an algebra over the Koszul resolution $\opd_{\infty}:=\Cobc \opd^{\textrm{!`}}$ of $\opd$ in the category of $\bialg$-operads. So, a structure of algebra over $\opd_{\infty}$ on an $\bialg$-module $A$ is equivalent to a morphism of $\bialg$-operads $\opd_{\infty}~:~=~\Cobc~\opd^{\textrm{!`}}~\rightarrow~\End_{A}$. Notice that a $\bialg$-$\opd_{\infty}$-algebra is a $\opd_{\infty}$-algebra endowed with a compatible $\bialg$-module structure.
\begin{Ex}
Any $\bialg$-$\opd$-algebra structure on a $\bialg$-module is a particular case of a $\bialg$-$\opd_{\infty}$-algebra structure.
\end{Ex}

\begin{pro}\label{inftystructure}
The set of $\bialg$-$\opd_{\infty}$-algebra structures is equivalently given by
$$\Hom_{\mathsf{dg \, }\bialg\textsf{-}\mathsf{Op}}(\Cobc \opd^{\textrm{!`}}, \End_{A}) \cong \Tw_{\bialg}(\opd^{\textrm{!`}},\End_{A}) \cong \Hom_{\mathsf{dg \, }\bialg\textsf{-}\mathsf{coOp}}(\opd^{\textrm{!`}}, \Barc \End_{A}) \cong \mathrm{Codiff}_{\bialg}(\opd^{\textrm{!`}}(A)) \ , $$
where $\mathrm{Codiff}_{\bialg}(\opd^{\textrm{!`}}(A))$ denotes the set of codifferentials on $\opd^{\textrm{!`}}(A)$, that are $\bialg$-equivariant coderivations on $\opd^{\textrm{!`}}(A)$ squaring to zero. 
\end{pro}

\begin{proo}
The first two bijections are given by the adjunction of Proposition \ref{adjunction}. The proof of Proposition \ref{Hopdalg} implies that 
$$ \Hom_{\SHMod}(\opd^{\textrm{!`}},\End_{A}) \cong \Hom_{\bialg}(\opd^{\textrm{!`}}(A),A) \ . $$
Moreover, \cite[Proposition 6.3.17]{LodayVallette09} extends to $\bialg$-modules. So, if $\coopd$ is a dg $\bialg$-cooperad, $V$ is an $\bialg$-module and $\alpha : \coopd(V) \rightarrow V$ is $\bialg$-equivariant, then the unique coderivation of the cofee $\coopd$-coalgebra $\coopd(V)$ which extends $\alpha$ is given by a sum of composites of $\bialg$-equivariant maps. Thus, we get the following isomorphism:
$$\Hom_{\bialg}(\opd^{\textrm{!`}}(A),A) \cong \mathrm{Coder}_{\bialg}(\opd^{\textrm{!`}}(A)) \ . $$
We conclude with \cite[Proposition 10.1.19]{LodayVallette09} stating that an element in $\Hom_{\SHMod}(\opd^{\textrm{!`}},\End_{V})$ is a solution to the Maurer-Cartan equation if and only if the associated coderivation on $\opd^{\textrm{!`}}(V)$ squares to zero, which does not depend on the $\bialg$-module structure.
\end{proo}

This result provides us with four equivalent definitions of $\bialg$-$\opd_{\infty}$-algebra structures. Thus, when dealing with this algebraic structure, we can make an ad hoc choice of one of those definitions.\\
\\
By extension of the classical definition of \cite[Section 10.2.2]{LodayVallette09}, an \textit{$\infty$-$\bialg$-morphism} of $\bialg$-$\opd_{\infty}$-algebras is morphism $A \rightsquigarrow B $ of dg $\bialg$-$\opd^{\textrm{!`}}$-coalgebras
$$F :(\opd^{\textrm{!`}}(A), d_{\varphi}) \rightarrow (\opd^{\textrm{!`}}(B), d_{\psi}) \ ,$$
where $\varphi$ and $\psi$ are the twisting $\bialg$-morphisms defining the structure of $\bialg$-$\opd_{\infty}$-algebra on $A$ and $B$ respectively. 
We denote by $\infty$-$\opd_{\infty}$-$\bialg$-$\mathsf{alg}$ the category of $\bialg$-$\opd_{\infty}$-algebras and their $\infty$-$\bialg$-morphisms, where the composite of two $\infty$-$\bialg$-morphisms is defined as the composite of the associated morphisms of dg $\bialg$-$\opd^{\textrm{!`}}$-coalgebras.\\
Equivalently, an $\infty$-$\bialg$-morphism of $\bialg$-$\opd_{\infty}$-algebras is an $\bialg$-equivariant $\infty$-morphism of the underlying $\opd_{\infty}$-algebras. A $\infty$-$\bialg$-morphism of $\bialg$-$\opd_{\infty}$-algebras is called an \textit{$\infty$-$\bialg$-isomorphism} (resp. \textit{$\infty$-$\bialg$-quasi-isomorphism}) if so is its first component $A \rightarrow B$.

\subsection{Homotopy transfer theorem}\label{HTT}

Let $(B,d_{B})$ be a homotopy retract of $(A,d_{A})$  in the category of $\bialg$-modules, that is a homotopy retract of chain complexes 
$$\xymatrix{  *{ \quad \ \  \quad (A, d_A)\ }  \ar@(dl,ul)[]^{h} \ar@<1ex>[r]^>>>>>{p} &  (B,d_{B})  \ar@<1ex>[l]^>>>>>{i}} \ ,$$
 with
$$\Id_{A} - ip = d_{A}h + hd_{A} , $$
such that the maps $i$, $p$ and $h$ are $H$-equivariant. We assume that $i$ is a quasi-isomorphism.

\begin{thm}\label{HTT}
Let $\opd$ be a Koszul $\bialg$-operad. Let $(B,d_{B})$ be a homotopy retract of $(A,d_{A})$  in the category of $\bialg$-modules. \\
Any $\bialg$-$\opd_{\infty}$-algebra structure on $A$, defined by generating operations $\lbrace m_\mu:A^{\otimes n}\rightarrow A, \mu \in \opd_{\infty}\rbrace$, can be transferred into a $\bialg$-$\opd_{\infty}$-algebra structure on $B$, which extends the  transferred operations $p \, m_{\mu}\,  i^{\otimes n}:B^{\otimes n} \rightarrow B$, and such that $i$ extends to an $\infty$-$\bialg$-quasi-isomorphism.
\end{thm}

\begin{proo}
We check that each step of the proof of Theorem 10.3.2 in \cite{LodayVallette09} extends to $\bialg$-modules. Since the maps $h$, $i$ and $p$ are $\bialg$-equivariant, then the map $ \Psi: \Barc\End_{A}\rightarrow \Barc\End_{B}$, defined in \cite[Section 10.3.3]{LodayVallette09}, is a morphism of $\bialg$-cooperads. To prove   that $\Psi$ is compatible with the differential structure, we do not care about the $\bialg$-module structure so it follows from \cite[Theorem 5.2]{VanDerLaan03}. Then, as a consequence of the Bar-Cobar adjunction, the composite $\opd^{\textrm{!`}} \rightarrow  \Barc\End_{A} \rightarrow \Barc \End_{B}$ defines an $\bialg$-$\opd_{\infty}$-algebra structure on $B$. Similarly, to be an $\infty$-$\bialg$-quasi-isomorphism does not depend on the $\bialg$-module structure therefore, by \cite[Theorem 10.3.11]{LodayVallette09}, we just have to prove that the map $i_{\infty}$, defined there, commutes with the $\bialg$-module structure. This is the case since it is a composite of such maps.
\end{proo}

In particular, the transferred structure and the $\infty$-$\bialg$-quasi-isomorphism are both given by the explicit tree-wise formulae of \cite[Section 10.3.10]{LodayVallette09} and of  \cite[Theorem 10.3.6]{LodayVallette09}. 

\begin{thm}
Let $\opd$ be a Koszul $\bialg$-operad. Let $(A,d_{A})$ be a  $\opd \rtimes \bialg$-algebra and $(B,d_{B})$ be a homotopy retract of $(A,d_{A})$ in the category of $\bialg$-modules. Then, $B$ inherits a structure of $H\text{-}\opd_{\infty}$-algebra, which extends the transferred operations, and such that $i$ extends to an $\infty$-quasi-isomorphism of $\bialg$-$\opd_{\infty}$-algebras.
\end{thm}

\begin{proo}
By Proposition \ref{semidirectprd}, the $\opd \rtimes \bialg$-algebra structure on $A$ is equivalent to an $\bialg$-$\opd$-algebra structure. Then, we apply Theorem \ref{HTT}.
\end{proo}

\begin{Rq} When the operad $\opd\rtimes \bialg$ is Koszul, one can apply the HTT for operads, but the $(\opd\rtimes\bialg)_{\infty}$-algebra structure can be very complex. If we can do the transfer in the category of $\bialg$-modules, then the transferred structure reduces to a much simpler structure. On the opposite, when the operad $\opd\rtimes \bialg$ is not Koszul, the classical method to transfer an algebraic structure must be improved. While, if the $\bialg$-operad $\opd$ is itself Koszul, then we can apply our version of the HTT, assuming that the homotopy retract is compatible with the $\bialg$-action.
\end{Rq}

\subsection{Homotopy theory for algebras over an $\bialg$-operad }\label{Hhoalg}

The results of \cite[Chapter 11]{LodayVallette09} extend \emph{mutatis mutandis} to the framework of $\bialg$-operads. The objects and the maps defined there are compatible with the $\bialg$-module structure. For the homological considerations, the results still holds because it does not depend on the $\bialg$-action. In particular, we obtain the following results :

\begin{thm}{(Rectification)}
Let $\opd$ be a Koszul $\bialg$-operad. Any $\bialg$-$\opd_{\infty}$-algebra is naturally $\infty$-$\bialg$-quasi-isomorphic to a dg $\bialg$-$\opd$-algebra.
\end{thm}

\begin{thm}{(Equivalence between homotopy categories)}\label{equivhomotopycat}
Let $\opd$ be a Koszul $\bialg$-operad. The homotopy category of dg $\bialg$-$\opd$-algebras is equivalent to the homotopy category of $\bialg$-$\opd$-algebras with their $\infty$-$\bialg$-morphisms.
\end{thm}

\begin{Rq}
In particular, we obtain a new description of the homotopy category of $\opd\rtimes\bialg$-algebras, using a simpler class of $\infty$-morphisms than the one of \cite[Theorem 11.4.12]{LodayVallette09}. 
\end{Rq}

An advantage of this result lies in the ``invertibility'' of the $\infty$-$\bialg$-quasi-isomorphisms. Indeed, any $\infty$-$\bialg$-quasi-isomorphism admits an $\infty$-$\bialg$-quasi-isomorphism in the opposite direction, as in \cite[Theorem 10.4.7]{LodayVallette09}, while it is not the case of $\bialg$-quasi-isomorphisms. This enables us to prove the following result.

\begin{thm}{($\bialg$-quasi-iso vs $\infty$-$\bialg$-quasi-iso)}\label{qivsinftyqi}
Let $\opd$ be a Koszul $\bialg$-operad and $A$, $B$ be two dg $\bialg$-$\opd$-algebras. The following assertions are equivalent:
\begin{enumerate}
\item there exists a zig-zag of $\bialg$-quasi-isomorphisms of dg $\bialg$-$\opd$-algebras 
$$ \xymatrix{A & \bullet \ar[l]_-{\sim} \ar[r]^-{\sim} & \bullet & \bullet \cdots  \bullet \ar[l]_-{\sim}   \ar[r]^-{\sim} & B} \ , $$
\item there exists an $\infty$-$\bialg$-quasi-isomorphism of dg $\bialg$-$\opd$-algebras $ A \overset{\sim}{\rightsquigarrow} B $.
\end{enumerate}
\end{thm}

\begin{Rq}
Since an $\bialg$-$\opd$-algebra structure is equivalent to a $\opd\rtimes\bialg$-algebra one, the theorem gives a new and simpler way to prove formality results for $\opd\rtimes\bialg$-algebras.
\end{Rq}

\subsection{Example}

%
%
%
%
%

When we encode BV-algebra structures with the operad $\BV$, a notion of homotopy BV-algebra structure is given in \cite{GalTonVal09}. We describe the homotopy structure obtained when we use the $\algnbduaux$-operad $\Gerst$ instead, and compare these two.\\


\begin{pro}\label{DGinfalg}
A homotopy $\algnbduaux$-$\Gerst$-algebra is a mixed chain complex $(A,d_{A},\unop_{A})$ endowed with a homotopy Gerstenhaber algebra structure such that 
$$ \unop_{A} \, m_{p_{1},\ldots,p_{t}}-(-1)^{t-2}\sum_{i=1}^{n}m_{p_{1},\ldots,p_{t}}\circ_{i} \unop_{A} = \displaystyle{\sum_{k=1}^{t}\sum_{p'_{k}+p''_{k}=p_{k}}} (-1)^{\varepsilon_{k}+1} m_{p_{1},\ldots,p'_{k},p''_{k},\ldots,p_{t}} \ , $$
where $n=p_{1}+\cdots+p_{t}$ and where $\unop_{A}:A\rightarrow A$ is the unary operator provided by the $\algnbduaux$-module structure on $A$.
\end{pro}

\begin{proo}
By Proposition \ref{inftystructure}, a structure of $\algnbduaux$-$\Gerst_{\infty}$-algebra on a mixed chain complex $A$ is given by a twisting $\algnbduaux$-morphism $\Gerst^{\textrm{!`}} \rightarrow \End_{A}$. It is a twisting morphism $\Gerst^{\textrm{!`}} \rightarrow \End_{A}$ between the underlying $\gsym$-modules, that is a homotopy Gerstenhaber structure as defined in~\cite[Proposition 16]{GalTonVal09}, which commutes with the action of $\algnbduaux$.    This condition is then a consequence of Proposition \ref{explicitDcoopstructure}, according to the study of homotopy Gerstenhaber algebras done in~\cite[Section 2.1]{GalTonVal09}.
\end{proo}

As conjectured at the end of \cite[Section 2.4]{GalTonVal09}, we have the following relation with the homotopy BV-algebras.

\begin{pro}\label{homBVtrivialmaps}
A homotopy BV-algebra structure (as defined in \cite[Theorem 20]{GalTonVal09}), such that all the operations are trivial except the maps $m^{0}_{p_{1},\ldots,{p_{t}}}$ and the map $m^{1}_{1}$, is equivalent to a $\algnbduaux$-$\Gerst_{\infty}$-algebra.
\end{pro}

\begin{proo}
Looking at the explicit description of the homotopy BV-algebra structure given in Theorem 20 of \cite{GalTonVal09}, all the relations become trivial except the relations $R^{0}_{p_{1},\ldots,p_{t}}$, defining the homotopy Gerstenhaber structure, the relations $R^{1}_{p_{1},\ldots,p_{t}}$, corresponding to those given in Corollary~\ref{DGinfalg}, and the relation $R^{2}_{1}$, which means that $m^{1}_{1}$ squares to zero.
\end{proo}

From now on, we call a $\algnbduaux$-$\Gerst_{\infty}$-algebra by a \emph{strict homotopy BV-algebra}, since a $\algnbduaux$-$\Gerst_{\infty}$-algebra is a homotopy BV-algebra such that the operator $\unop$ strictly satisfies the relations of a BV operator. 

\begin{Rq}
A strict homotopy BV-algebra is a homotopy BV-algebra such that any elements of $\BV^{\textrm{!`}}$ (as defined in \cite[Section 1]{GalTonVal09}) containing a vertex decorated by $\unop$, except $\unop$ itself, acts as zero. If we ask, moreover, that $\unop$ acts as zero, we get the notion of \emph{strongly trivialized homotopy BV-algebras}, defined in \cite{DrummondCole12}. In this article, the author proved that any strongly trivialized homotopy BV structure on a chain complex induces a structure of homotopy hypercommutative algebra on that complex. Thus, we get another description of the homotopy hypercommutative structure. 
\end{Rq}

\begin{thm}\label{strictBVinf}
Let $(B,d_{B},\unop_{B})$   be a homotopy retract of $(A,d_{A}, \unop_{A})$ in the category of mixed chain complexes. If $A$ is endowed with a BV-algebra structure (or with a strict homotopy BV-algebra structure), then $B$ inherits a structure of strict homotopy BV-algebra which extends the naive transferred operations, and such that $i$ extends to an $\infty$-$\algnbduaux$-quasi-isomorphism of strict homotopy BV-algebras.
\end{thm}

\begin{proo}
By Proposition \ref{BValg}, $A$ is a $\algnbduaux$-$\Gerst$-algebra. Then we apply Theorem \ref{HTT}.
\end{proo}

When the  homotopy retract lives in the category of mixed chain complexes, this theorem shows that the 
transferred homotopy $BV$-algebra structure reduces to a strict homotopy $BV$-algebra, i.e a homotopy 
BV-algebra structure without higher homotopies arising from $\unop$ and its relations with the product and with the bracket.


\begin{pro}\label{BVside}
Let $(B,d_{B})$ be a homotopy retract of $(A,d_{A})$, in the category of chain complexes, which satisfies the following \emph{side conditions}:
$$ h^{2}=hi=ph=0 \ . $$
If $(A,d_{A}, \prd{}{},\bracket{\,}{},\unop)$ is a BV-algebra structure on $A$ such that 
$$[\unop , h ] =0\ ,$$
then $B$ inherits a strict homotopy BV-algebra structure, induced by the BV-algebra structure on $A$, and such that $i$ extends to an $\infty$-$\algnbduaux$-quasi-isomorphism of strict homotopy BV-algebras. 
\end{pro}

\begin{proo}
We endow the chain complex $(B,d_{B})$ with the unary operator $\widetilde{\unop{}}$ induced by $\unop$, and defined by $\widetilde{\unop{}}:= p\unop i $. Then, the side conditions imply that we have actually a homotopy retract between $(A,d_{A}, \unop)$ and $(B,d_{B}, \widetilde{\unop{}})$ in the category of mixed chain complexes. We conclude with Theorem \ref{strictBVinf}.
\end{proo}

This is the result which answers the main question raised at the begining of this paper.


When $B$ is the homology groups $\mathrm{H}^{\bullet}(A,d)$ of the complex $(A,d)$, it is always possible to build a homotopy retract which satisfies the side conditions: it is called the \emph{Hodge decomposition}, see \cite[Lemma 9.4.7]{LodayVallette09}. More precisely, the chain complex $(A,d)$ splits into a direct sum of graded spaces as follows: 
$$ A = H \oplus B \oplus C \ , $$
where $B= \mathrm{Ker}(d)\cap \mathrm{Im}(d)$, $H\oplus B = \mathrm{Ker}(d)$ and $H\cong \mathrm{H}^{\bullet}(A,d)$. Then, the homotopy retract is defined as follows : 
\begin{itemize}
\item[$\bullet$] $i$ is the inclusion of $H$ in $A$,
\item[$\bullet$] $p$ is the projection of $A$ onto $H$,
\item[$\bullet$] and $h$ is equal to $0$ on $H\oplus C$ and is equal to $d^{-1}$ on $B$.
\end{itemize}  
In this case, the homotopy retract satisfies the aforementioned side conditions.

%

\begin{Ex}
We consider the following non-unital dg commutative algebra $A$ generated by 
$$ x_{3}, y_{3},t_{3}, \xi_{4}, \omega_{5}, z_{7}, u_{7} \text{ and } v_{8} \ , $$
where the subscript gives the homological degree, such that the product by $u$, by $v$, by $\xi$ and by $\omega$ is equal to zero. The algebra $A$ is finite dimensional and spanned by the elements $x$, $y$, $t$, $\xi$, $\omega$, $xy$, $xt$, $yt$, $u$, $z$, $v$, $xz$, $yz$, $tz$, $xyz$, $xyt$ and $xyzt$. We summarize in the following picture the definition of the differential and of the BV operator on the generators :
$$ \xymatrix@R=7pt@M=5pt@C=15pt{ 0 & 3 & 4 & 5 & 6 & 7 & 8 & 9 & 10 & 13 & 16 \\ & x & \xi \ar@{|->}[r]^-{\unop{}} & \omega & xy \ar@{|->}[r]^-{\unop{}} & u & v \ar@{|->}[l]_{d} & xyt & xz & xyz & xyzt \\  & y & & &  xt & z \ar@{|->}[lu]^{d} \ar@{|->}[ru]_-{-\unop{}}  & & & yz & \\ & t & & & yt & & & & zt } \ , $$
where an element is sent to $0$ if nothing else is specified.
\end{Ex}

\begin{pro}
The aforementioned algebra $(A,d,\unop)$ is a dg BV-algebra. Moreover, its homology $H^{\bullet}(A,d)$ is endowed with a strict homotopy BV-algebra structure, such that the inclusion $H^{\bullet}(A,d) \hookrightarrow A$ extends to an $\infty$-quasi-isomorphism of strict homotopy BV-algebras.
\end{pro}

\begin{proo}
It is straightforward to prove that $A$ together with $d$ and $\unop$ forms a dg BV-algebra. The Hodge decomposition of $A$ is given by 
$$ \xymatrix@R=7pt@M=5pt@C=15pt{ n & 3 & 4 & 5 & 6 & 7 & 8 & 9 & 10 & 13 & 16 \\ H_{n} & x,y,t & \xi \ar@{|->}[r]^-{\unop{}} & \omega & xt, yt &  &  xyt & & xz, yz, tz & xyz & xyzt \\  B_{n} & & & & xy \ar@{|->}[r]^-{\unop{}} \ar@{|->}[rd]_-{h=d^{-1}}& u \ar@{|->}[rd]^{h=d^{-1}}  & & & & \\ C_{n} & & & & & z \ar@{|->}[r]_-{-\unop{}} & v & &  } \ , $$
where an element is sent to $0$ if nothing else is specified. The contracting homotopy and $\unop$ commute. Then, we apply Proposition \ref{BVside}. 
\end{proo}
Note that the product induced in homology is associative, since the differential on $H^{\bullet}(A,d)$ is $0$, but the first homotopy $m^{0}_{3}$ is not equal to $0$. Moreover, the BV operator does not vanish in homology.

\begin{Rq}
Introduced in \cite{DGMS75} to study the differential forms of compact K\"ahler manifolds and used in \cite{BarKon98} for the Dolbeault complex of Calabi-Yau manifolds, the $dd^{c}$-lemma implies the condition $[\unop , h ] =0$, but is a strong condition in our setting. Indeed, if the $dd^{c}$-lemma is satisfied, then the operator $\unop$ vanishes on the homology groups.
\end{Rq}

Encoding the BV-algebra structure in the $\algnbduaux$-operad $\Gerst$ allows us to deal with a notion of homotopy BV-algebra which is simpler than the one given in \cite{GalTonVal09}. The associated notion of $\infty$-morphism is also simplified.

\begin{pro}
An $\infty$-$\algnbduaux$-morphism of (strict homotopy) BV-algebras is an $\infty$-morphism of the underlying (homotopy) Gerstenhaber algebras which commutes to the extra action of the unary operator, provided by the (strict homotopy) BV-algebra structures.
\end{pro}

\begin{proo}
It is a consequence of the definition of an $\infty$-morphism in the category of $\algnbduaux$-$\Gerst_{\infty}$-algebras.
\end{proo}

\begin{cor}\label{zigzagBV}
The homotopy category of BV-algebras is equivalent to the homotopy category of strict homotopy BV-algebras with their $\infty$-morphisms. Moreover, for any BV-algebras $A$ and $B$, the following assertions are equivalent:
\begin{enumerate}
\item there exists a zig-zag of quasi-isomorphisms of BV-algebras 
$$ \xymatrix{A & \bullet \ar[l]_-{\sim} \ar[r]^-{\sim} & \bullet & \bullet \cdots  \bullet \ar[l]_-{\sim}   \ar[r]^-{\sim} & B} \ , $$
\item there exists an $\infty$-quasi-isomorphism of Gerstenhaber algebras $ A \overset{\sim}{\rightsquigarrow} B $, commuting with the unary operator provided by the BV-algebra structures.
\end{enumerate}
\end{cor}

\begin{proo}
Since the operad $\Gerst$ is Koszul, we can apply Theorem \ref{equivhomotopycat} and Theorem \ref{qivsinftyqi}, combined with the previous Proposition. 
\end{proo}

This theorem provides a new way to prove the existence of a chain of quasi-isomorphisms between BV-algebras. In particular, it could help to prove the formality of a BV-algebra. Moreover, we could use it to study the interpretation of the mirror symmetry conjecture of Kontsevich \cite{Kontsevich95} in terms of BV-algebras.

The mirror symmetry conjecture says that a certain A-type theory of a Calabi-Yau manifold can be identified with a B-type theory on a mirror manifold. In \cite{CaoZhou01}, the authors interpret this conjecture as the existence of a quasi-isomorphism of dg BV-algebras from the De Rham complex of a Calabi-Yau manifold $\mathcal{M}$ to the Dolbeault complex of a dual Calabi-Yau manifold $\widetilde{\mathcal{M}}$:
$$ \left(\Omega^{n-\bullet}(\mathcal{M}), d_{DR},\wedge,\unop, \langle \ ; \, \rangle \right) \xrightarrow{\sim} \left( \Gamma ( \widetilde{\mathcal{M}},\wedge^{\bullet}\overline{T}^{\ast}_{\widetilde{\mathcal{M}}} \otimes \wedge^{\bullet}T_{\widetilde{\mathcal{M}}}),\overline{\partial}, \wedge, \mathrm{div}, \langle \ ; \, \rangle_{S} \right) \ . $$
Using Corollary \ref{zigzagBV}, it is equivalent to the existence of an $\infty$-quasi-isomorphism of the underlying Gerstenhaber algebras, which commutes strictly with the action of the unary operators $\unop$ and $\mathrm{div}$.

\DeactivateToc
\section*{Acknowlegments}
\ActivateToc

This paper constitutes a part of the author's Ph.D. thesis. I would like to thank my advisor, Bruno Vallette, for being so receptive and involved. I am grateful for our instructive discussions and for his fruitful ideas. I also would like to thank Joan Mill\`es for his reactivity to my many questions. The author is supported by the ANR grant HOGT, referenced as ANR-11-BS01-0002.

\bibliographystyle{plain}
\bibliography{bib}

\end{document}